\documentclass[11pt]{article}

\usepackage{amsmath}
\usepackage{amsfonts}
\usepackage{amssymb}
\usepackage{amsthm}
\usepackage{graphics}
\usepackage{amscd}
\usepackage{graphicx,epsfig}
\usepackage{color}
\usepackage[all]{xy}
\usepackage{url}

\DeclareMathOperator{\Div}{div}

\DeclareMathOperator{\argth}{argth}

\renewcommand{\epsilon}{\varepsilon}

\newcommand{\boD}{\mathcal{D}}

\newcommand{\R}{\mathbb{R}}

\renewcommand{\H}{\mathbb{H}}
\renewcommand{\S}{\mathbb{S}}
\newcommand{\N}{\mathbb{N}}

\newcommand{\eps}{\varepsilon}

\newcommand{\dd}{\mathrm{d}}

\newcommand{\Ome}{\Omega}

\newtheorem{defn}{Definition}
\newtheorem{thm}{Theorem}

\newtheorem{prop}[thm]{Proposition}
\newtheorem{lem}[thm]{Lemma}

\renewcommand{\phi}{\varphi}
\newcommand{\dis}{\displaystyle}

\newtheorem*{thm*}{Theorem}

\newtheorem{claim}[thm]{Claim}
\newtheorem*{claim*}{Claim}

\theoremstyle{remark}

\newtheorem*{rem*}{Remark}

\newcounter{remark}

\newcounter{case}

\newcounter{construction}

\newcounter{fact}

\newcounter{step}

\newcommand{\bnabla}{\overline{\nabla}}

\title{Cylindrically bounded constant mean curvature surfaces in $\H^2\times\R$}
\author{Laurent Mazet\thanks{The author was partially supported by the
ANR-11-IS01-0002 grant.}}
\date{}

\begin{document}

\maketitle

\begin{abstract}
In this paper we prove that a properly embedded constant mean curvature surface
in $\H^2\times\R$ which has finite topology and stays at a finite distance from
a vertical geodesic line is invariant by rotation around a vertical geodesic
line.
\end{abstract}

\section{Introduction}
In 1988, W.~Meeks \cite{Mee} proved that, in $\R^3$, a properly embedded annulus
with non vanishing constant mean curvature (cmc, in the following) must stay at
a bounded distance from a straight line. Then N.~Korevaar, R.~Kusner and
B.~Solomon \cite{KoKuSo} proved that any properly embedded constant mean
curvature surface staying at a bounded distance of a straight line is
rotationally invariant and then a Delaunay surface.

These results imply that any end of a finite topology properly embedded cmc
surface is asymptotic to a Delaunay surface. So this allows a description of the
space of all properly embedded cmc surfaces with finite topology (see the paper
of R.~Kusner, R.~Mazzeo and D.~Pollack \cite{KuMaPo}).

N.~Korevaar, R.~Kusner and B.~Solomon used the Alexandrov reflection procedure to
prove their results. In fact, their proof works in higher dimension and also in
$\H^n$ \cite{KoKuMeSo}.

Recently the theory of cmc surfaces is developed in 3-dimensional homogeneous
spaces. One interesting case is the ambient space $\H^2\times\R$. The group of
isometries of $\H^2\times\R$ possesses rotations around vertical geodesic lines
$p\times\R$ where $p$ is a point of $\H^2$. So one can look for embedded cmc
surfaces invariant by rotation around such a geodesic line, this was done by
Hsiang and Hsiang in \cite{HsHs}. For any mean curvature $H_0$, they find
rotationally invariant surfaces but these surfaces stay at a bounded distance
from the vertical geodesic axis only for $H_0>1/2$. In the sequel, we will focus
on the case $H_0>1/2$. Among these rotationally invariant surfaces, certain are
spheres and, actually, they are the only compact embedded cmc surface in
$\H^2\times\R$ \cite{HsHs}. The other ones are periodic with respect to a vertical
translation. For example, we have the vertical cylinder $C\times\R$ where $C$ is
a circle in $\H^2$. So these surfaces correspond to the Delaunay surfaces in
$\R^3$ and they are also called Delaunay surfaces in $\H^2\times\R$.

A subset of $\H^2\times\R$ which is at a bounded distance from a vertical
geodesic line will be called \emph{cylindrically bounded}. So, as in the paper
of Korevaar, Kusner and Solomon: can we classify properly embedded cmc surface
that are cylindrically bounded? The main result of the paper
(Theorem~\ref{th:main}) gives an answer to this question.

\begin{quote}
\textbf{Theorem.} Let $\Sigma$ be a properly embedded cmc surface in
$\H^2\times\R$. If $\Sigma$ has finite topology and is cylindrically bounded,
$\Sigma$ is a Delaunay surface.
\end{quote}

So any cylindrically bounded cmc annulus in $\H^2\times\R$ is a Delaunay surface.
Let us make a remark about the limit case $H_0=1/2$. In $\R^3$, the results of
Meeks and Korevaar, Kusner and Solomon imply that any properly embedded cmc
$H_0>0$ annulus is rotational. For $H_0=0$, we also know that a properly
embedded minimal annulus is a catenoid and thus rotational (see Collin
\cite{Col2}). In $\H^2\times\R$, we do not have such a rigidity. Actually,
Cartier and Hauswirth \cite{CaHa} have recently constructed a properly embedded
cmc $H_0=1/2$ annulus in $\H^2\times\R$ whose ends are asymptotic to rotational
ones but which is not rotationally invariant.

We notice that a vertical cylinder $C\times\R$ has mean curvature larger than
$1/2$. Thus, any properly embedded cmc surface which is cylindrically bounded
has mean curvature $H_0>1/2$. This is a consequence of the half space theorem in
$\H^2\times\R$ (see \cite{NeRo3} and \cite{Maz15}, for example). So we can focus
ourselves only on the $H_0>1/2$ case.

The proof of the theorem is also based on the Alexandrov reflection technique
but the space of planar symmetries in $\H^2\times \R$ is smaller than in $\R^3$
(in $\R^3$ it is a $3$ dimensional space and in $\H^2\times\R$ it is only $2$
dimensional). So the ideas of Korevaar, Kusner and Solomon can not be applied.

First we remark that we already know that a compact embedded cmc surface is a
rotational sphere so we only consider non compact surfaces. For a non compact
cylindrically bounded cmc surface $\Sigma$ and a foliation of $\H^2\times\R$ by
vertical planes, we define a function $\alpha$ on $\R$ called the
\emph{Alexandrov function}. One proprety of this function is that, if $\alpha$
admits a maximum, 
$\Sigma$ is symmetric with respect to a vertical plane of the foliation. To
prove that $\Sigma$ is rotationally invariant, it suffices to prove that it
is symmetric with respect to a lot of vertical planes. So we want to prove that
any Alexandrov function has a maximum.

If the surface $\Sigma$ has finite topology, we know by previous results
\cite{HoLiRo} that it
has bounded curvature. So we can control the asymptotic behavior of the surface.
This gives informations about the behavior of the Alexandrov function $\alpha$ near
$\pm\infty$. Then we prove that $\alpha$ is decreasing near $+\infty$
and increasing near $-\infty$: this implies that $\alpha$ has a maximum. To
prove this monotonicity result we use a flux argument similar to the one of the
Positive flux lemma proved by Korevaar, Kusner, Meeks and Solomon in
\cite{KoKuMeSo}.

We notice that the proofs given in this paper work also in $\S^2\times\R$: they 
prove that a properly embedded constant mean curvature surface with finite
topology in $D\times\R$ is rotationally invariant where $D$ is geodesic disk in
$\S^2$ of radius less than $\pi/2$.

The paper is divided as follows. In Section~\ref{sec:preliminar}, we recall
several result concerning cylindrically bounded cmc surfaces in $\H^2\times\R$.
In Section~\ref{sec:delaunay}, we recall the construction of the Delaunay
surfaces. In Section~\ref{sec:alexandrov}, we define the Alexandrov function and
give its properties, we study the asymptotic behavior of an annular end of a
cylindrically bounded cmc surface. Finally we state our main result and give the
main steps of the proof. Section~\ref{sec:horigraph} is devoted to the study of
horizontal Killing graph, \textit{i.e.} cmc surfaces that are transverse to the
horizontal Killing vector field generating horizontal translations in
$\H^2\times\R$. The last section is devoted to the study of the monotonicity of
the Alexandrov function near $\pm\infty$. The main idea is to compare the flux
along the surface with the flux along a comparison surface which is constructed
as a horizontal Killing graph.

All along the paper, we use $z$ to denote the real coordinate in
$\H^2\times\R$. On $\H^2$, we consider the polar coordinates $(\rho,\theta)$
in Section~\ref{sec:delaunay} and a $(s,r)$ coordinate system in
Sections~\ref{sec:horigraph} and \ref{sec:monotonous}. This last coordinate
system is adapted to the Killing vector field that generates translations along
a geodesic line.

\section{Previous results}
\label{sec:preliminar}
In this section, we recall different results concerning properly embedded cmc
surfaces in $\H^2\times\R$. Most of them can be found in the papers of
D.~Hoffman, J.~de~Lira and H.~Rosenberg \cite{HoLiRo,Rog2} and the one by
N.~Korevaar, R.~Kusner and B.~Solomon \cite{KoKuSo}. We also explain the
convergence we will consider for sequences of cmc surfaces with bounded
curvature.


\subsection{Flux} 
This first subsection is devoted to the notion of flux for cmc surfaces that was
introduced in several preceding papers (see for example \cite{KoKuSo,HoLiRo}). 

\begin{figure}[h]
\begin{center}
\resizebox{0.6\linewidth}{!}{\input{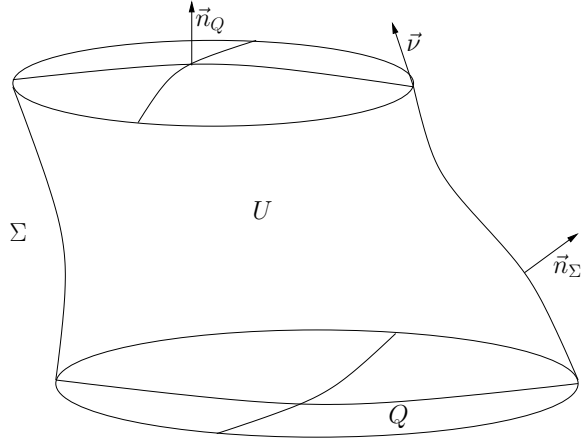}}
\caption{Computation of the flux}\label{fig:flux}
\end{center}
\end{figure}

Let $U$ be a bounded domain in $\H^2\times\R$ whose boundary
is the union of a smooth surface $\Sigma$ and a smooth surface $Q$ with
common boundary $\partial \Sigma=\partial Q$. So the boundary $\partial U$ is
piecewise smooth. Let us denote by $\vec n$ the unit outgoing normal along
$\partial U$ and $\vec n_\Sigma$ and $\vec n_Q$ the respective restriction of
$\vec n$ along $\Sigma$ and $S$ (see Figure~\ref{fig:flux}).

Let $Y$ be a Killing vector field of $\H^2\times\R$, we have 
$$
0=\int_U \Div_{\H^2\times\R} Y=\int_\Sigma Y\cdot \vec n_\Sigma+ \int_Q Y\cdot
\vec n_Q
$$
where $\Div_{\H^2\times\R}$ is the divergence operator on $\H^2\times\R$.

Let $\vec\nu$ be the outgoing unit conormal to $\Sigma$ along $\partial \Sigma$.
Along $\Sigma$ the vector field $Y$ can be decomposed into the sum of a tangent
part $Y^\top$ and a normal part $Y^\bot$. We have
$$
0=\int_\Sigma \Div_\Sigma Y=\int_\Sigma \Div_\Sigma Y^\top+ \int_\Sigma \Div_\Sigma
Y^\bot = \int_{\partial \Sigma} Y\cdot \vec\nu +\int_\Sigma 2HY \cdot \vec n_\Sigma
$$
where $\Div_\Sigma$ is the divergence operator on $\Sigma$ acting on any vector
field of $\H^2\times\R$ defined along $\Sigma$ and $H$ is the mean curvature of
$\Sigma$ computed with respect to $-\vec n_\Sigma$.

Thus, if $\Sigma$ has constant mean curvature $H_0$, we have 
\begin{equation}\label{eq:defflux}
0=\int_{\partial \Sigma} Y\cdot \vec\nu- 2H_0 \int_Q Y\cdot \vec n_Q.
\end{equation}

Now let us consider $\Sigma$ a cmc $H_0$ surface and $\gamma$ a smooth closed
curve in $\Sigma$. Let $Q$ be a smooth surface in $\H^2\times\R$ with boundary
$\gamma$. For a Killing vector field $Y$, we define the quantity 
$$
F_\gamma(Y)=\int_\gamma Y\cdot \vec\nu -2H_0\int_Q Y\cdot \vec n_Q
$$
where, as above, $\vec\nu$ is the conormal unit vector field of $\Sigma$ along
$\gamma$ and $\vec n_Q$ is a unit normal along $Q$ (the choice of these normal
vectors is consistent with the above computations). Because of Formula
\eqref{eq:defflux},
the quantity $F_\gamma(Y)$ does not depend on the choice of $Q$ and depends only
on the homology class of $\gamma$ in $\Sigma$. $F_\gamma(Y)$ is called the flux
of $\Sigma$ along $\gamma$ in the direction $Y$.

In fact, the map $Y\mapsto F_\gamma(Y)$ is linear so it can be seen as an element
of the dual of the vector space of Killing vector fields. This element is called
the flux of $\Sigma$ along $\gamma$.


\subsection{Linear area growth}

In this subsection, we recall a result concerning the area growth of a
cylindrically bounded properly embedded cmc surface in $\H^2\times \R$. 

Let $\Sigma$ be a properly embedded cmc surface with possibly non-empty compact
boundary. We say that $\Sigma$ has linear area growth if there exists two
constants $\alpha$ and $\beta$ such that $Area(\Sigma \cap \{ a\le z\le b\})\le
\alpha(b-a)+\beta$ for any $a\le b\in \R$.

We then have the following result.
\begin{prop}[Corollary 1 in \cite{HoLiRo} or Theorem 14 in
\cite{Rog2}]\label{prop:linarea}
Let $\Sigma\subset\H^2\times\R$ be a properly embedded cmc surface with possibly
non-empty compact boundary. If $\Sigma$ is cylindrically bounded, $\Sigma$ has
linear area growth.
\end{prop}

We notice that in this result we do not make any hypothesis on the topology of
the surface $\Sigma$. Besides, the original proof in \cite{HoLiRo} contains a
mistake, it has been corrected by Rosenberg in \cite{Rog2}.

\subsection{The height function}

On a properly embedded surface $\Sigma$ in $\H^2\times\R$ the restriction to the
surface of the real coordinate $z$ is called the height function on $\Sigma$. For a
cmc $H_0$ surface, this height function has several properties.

\begin{lem}\label{lem:heightbound}
Let $\Sigma$ be a non compact properly embedded cylindrically bounded cmc $H_0$
surface in $\H^2\times\R$. The height function can not be neither bounded from
above nor bounded from below.
\end{lem}

The idea of the proof can be found in the proof of Theorem~1.1 in \cite{NeRo2}
or Proposition~2 in \cite{HoLiRo}.

\begin{proof}
For example, let us assume that the height function is bounded from below. We
can apply Alexandrov reflection technique with respect to horizontal slice
$\H^2\times\{t\}$. Since $\Sigma$ is bounded from below, for $t$ small, $\Sigma$
does not intersect $\H^2\times\{t\}$. Thus we can start the reflection procedure
up to a first contact point. But since $\Sigma$ is non compact, a first contact
point can not exist. Indeed, if there is a first contact point, the maximum
principle would imply that $\Sigma$ is symmetric with respect to some
$\H^2\times\{t_0\}$ and thus compact. 

So the reflection procedure can be done
for any $t$ and this implies
that the part of $\Sigma$ below $\{z=t\}$ is a vertical graph with boundary in
$\{z=t\}$. The height of such a vertical graph is bounded from above by a
constant that depends only on $H_0$ (see \cite{AlEsGa}); thus the Alexandrov
procedure has to stop. We get a contradiction.
\end{proof}

When the surface $\Sigma$ is an annulus we have the following property (see
Lemma~4.1 in \cite{KoKuSo}).

\begin{lem}\label{lem:compannulus}
Let $H_0>1/2$ be a real number, there exists an $M>0$ that depends only on
$H_0$ such the following is true. Let $A$ be an embedded annulus $A\subset
\H^2\times\R$ of constant mean curvature $H_0$ with boundary outside
$\H^2\times[0,M]$. Then $A\cap \H^2\times[0,M]$ has at most one connected
component $\tilde A$ such that $z(\tilde A)=[0,M]$.
\end{lem}


\subsection{Uniform curvature estimate}

In this subsection, we recall an estimate of the norm of the second
fundamental form of a cylindrically bounded properly embedded cmc surface in
$\H^2\times \R$. Precisely, Hoffman, de~Lira and Rosenberg proved the following
result.
\begin{prop}[Theorem 3 in \cite{HoLiRo}]\label{prop:boundcurv}
Let $\Sigma$ be a properly embedded cmc surface with finite topology and
possibly non-empty compact boundary. If $\Sigma$ is cylindrically bounded, the
norm of the second fundamental form $|A|$ is bounded on $\Sigma$.
\end{prop}


\subsection{Convergence of sequences of cmc surfaces}\label{sec:conv}

In the following, we will consider sequences of cmc surfaces coming from the
translations of a given cmc surfaces with bounded curvature. In this subsection,
we explain how these sequences converge to a limit cmc surface. The surfaces we
consider have a uniform curvature bound.

Considering the normal coordinates around a point in $\H^2\times\R$, a
constant mean curvature surface in
$\H^2\times\R$ can be viewed has an immersed surface in $\R^3$. A bound of its
second fundamental form in $\H^2\times\R$ is equivalent to a bound in $\R^3$. We
have the following classical result.

\begin{prop}\cite{KoKuSo,PeRo}\label{prop:localbehav}
Let $\Sigma$ be an immersed surface in $\R^3$ whose second fundamental form
satisfies $|A|\le 1/(4\delta)$ for some $\delta>0$. Then for any $x\in\Sigma$
with $d(x,\partial\Sigma)> 4\delta$ there is a neighborhood of $x$ in $\Sigma$
which is a graph of a function $u$ over the Euclidean disk of radius $\sqrt
2\delta$ centered at $x$ in the tangent plane to $\Sigma$ at $x$. Moreover 
$$
|u|<2\delta,\ \ |\nabla u|<1,\text{ and } |\nabla^2 u|< \frac 1\delta
$$
\end{prop}

Let $\Sigma$ be a cmc $H_0$ surface ($H_0\neq 0$) that bounds a domain $D$ in
$\H^2\times\R$ ($D$ is in the mean convex side of $\Sigma$) and with a uniform
curvature bound. Let $p\in\Sigma$ be a point and $\gamma$ be the geodesic line
starting from $p$ in the direction of the mean curvature vector. Let $q$ be 
the first point where $\gamma$ meets $\Sigma$ (if it exists). If $q$ is close to
$p$, Proposition~\ref{prop:localbehav} implies that $\gamma$ is close to be
normal to $\Sigma$ at $q$. Moreover, since $\gamma$ is in $D$ between $p$ and $q$, the
mean curvature vector to $\Sigma$ at $q$ points in the opposite direction to the
velocity vector of $\gamma$ at $q$. But since $H_0>0$, this implies that $p$ and
$q$ can not be too close from each other. More precisely, we have the following
result that gives a local description of the surface.

\begin{prop}\label{prop:localbehav2}
Let $\Sigma$ be an embedded cmc $H_0$ surface in $\H^2\times\R$ ($H_0>0$) which
bounds a domain $D$ in its mean convex side. Moreover we assume that its second
fundamental form satisfies $|A|\le k$ for some $k>0$. Then there exists
$R=R(k)>0$ such that the following is true. For any $p\in \H^2\times\R$, there
exists at most two topological disks $\Delta^1$ and $\Delta^2$ in $\Sigma\cap
B(p,2R)$ such
that $\Sigma\cap B(p,R)=(\Delta^1\cup\Delta^2)\cap B(p,R)$. Moreover, when there
is two disks, the domain between the two disks in $B(p,R)$ is outside $D$
\end{prop}

\begin{figure}[h]
\begin{center}
\resizebox{0.8\linewidth}{!}{\input{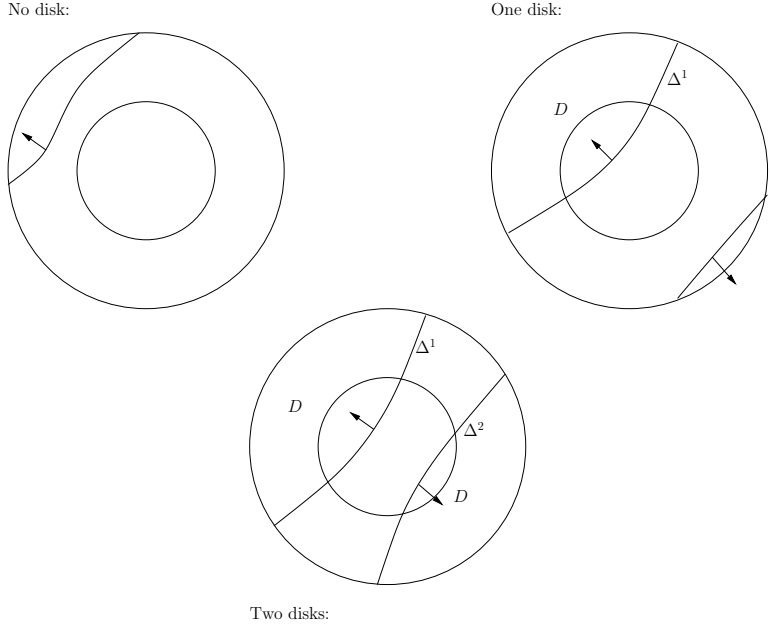}}
\caption{The local description of a cmc surface}\label{fig:descrip}
\end{center}
\end{figure}

The above proposition says that a local description of the surface $\Sigma$ is
given by one of the pictures in Figure~\ref{fig:descrip}.

Now let us assume that we have a sequence $(\Sigma_n)$ of properly embedded cmc
$H_0$ surfaces with uniformly bounded curvature and bounding a domain $D_n$ in
their mean convex side. For each compact subset $K$ of $\H^2\times \R$, $(K\cap
\Sigma_n)$ is a sequence of compacts sets so by a diagonal process we can assume
that the sequence $(\Sigma_n)$ converges in the Hausdorff topology restricted to
each $K$ to a limit that we
call $\Sigma$. Let $p$ be in $\Sigma$, there exists $p_n\in\Sigma_n$ such that
$p_n\rightarrow p$. From Proposition~\ref{prop:localbehav}, $p_n$ is the
center of a geometrically controlled disk $\Delta_n^1$. By considering a
subsequence we can assume that $\Delta_{n'}^1$ converges to a cmc $H_0$ disk
$\Delta^1$. We then have $\Delta^1\subset\Sigma$. 

Let $R$ be given by Proposition~\ref{prop:localbehav2}. If $\Sigma\cap
B(p,R)=\Delta^1\cap B(p,R)$, this implies that the whole sequence $\Delta_n^1$
converges to $\Delta^1$. If there is $q\in (\Sigma\cap B(p,R))\setminus
\Delta^1$, there exist $q_n\in \Sigma_n$ and $q_n$ is in the second disk
$\Delta_n^2$ given by Proposition~\ref{prop:localbehav2}. By considering a
subsequence we can assume $\Delta_{n'}^2$ converges to a cmc $H_0$ disk
$\Delta^2$. We notice that $\Delta^2$ is different from $\Delta^1$ since it
contains $q$. We then have $\Delta^2\subset \Sigma$. From
Proposition~\ref{prop:localbehav2}, we can be sure that $\Sigma\cap
B(p,R)=(\Delta^1\cup \Delta^2)\cap B(p,R)$. In fact this implies that the whole
sequence $(\Delta_n^2)$ converges to $\Delta^2$.

To complete our local description of $\Sigma$, we notice that $\Delta^1$ and
$\Delta^2$ can touch each other but, at these contact points, the mean curvature
vector have opposite value. Besides $\Sigma$ bounds the domain $D$ which is
constructed has the limit of the domain $D_n$. This prove that $\Sigma$ is what
we call a weakly embedded cmc $H_0$ surface.

\begin{defn}
A properly immersed cmc surface $\Sigma$ is said to be weakly embedded if there
exists an open set $\Ome$ such that $\Sigma$ is the boundary of $\Ome$ and the
mean curvature vector of $\Sigma$ points into $\Ome$. 
\end{defn}

As an example, two rotational cmc spheres that are tangent form a weakly
embedded surface. The union of two vertical cylinders tangent along a common
vertical geodesic line is also weakly embedded.

We say that a weakly embedded surface is connected if the underlying abstract
surface $\Sigma$ is connected. As an example, two rotational cmc sphere tangent
at a point is not a connected weakly embedded surface.

Finally, we remark that we have proved that the sequence $\Sigma_n$ converges
smoothly in any compact sets to the surface $\Sigma$.

\section{Delaunay surfaces}
\label{sec:delaunay}
In this section, we briefly recall the construction of the embedded cmc surfaces
that are rotationally invariant around a vertical axis. We only focus on
surfaces with $H>1/2$ (see \cite{HsHs} and \cite{PeRi} for more details). 

\begin{figure}[h]
\begin{center}
\resizebox{0.8\linewidth}{!}{\input{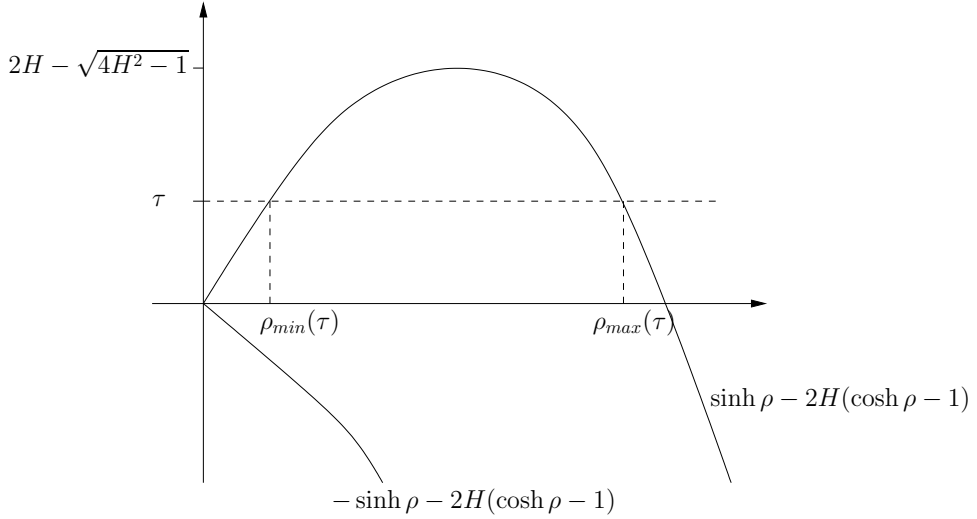}}
\caption{The domain of definition of solutions}\label{fig:taurho}
\end{center}
\end{figure}

Let $(\rho,\theta)$ be the polar coordinates on $\H^2$ so the metric is
$\dd \rho^2+ (\sinh \rho)^2 \dd\theta^2$. We look for surfaces of revolution; so
we can look for the graph of a function $u=f(\rho)$, we will orient this graph
using the upward pointing unit normal. Since we want a constant mean
curvature $H$ graph, the function $f$ satisfies the following equation
$$
\dis\frac{f'}{\sqrt{1+f'^2}}\sinh \rho-2H(\cosh\rho-1)=\tau
$$
where $\tau$ is a constant. In order to have a solution, $\tau$ has to satisfy
$$
-\sinh\rho-2H(\cosh \rho-1)\le \tau\le \sinh \rho-2H(\cosh \rho-1)
$$
The graphs of the functions in the left-hand and right-hand sides are given in
Figure~\ref{fig:taurho}. So in order to have a
solution with non empty definition set, $\tau$ has to be chosen less than
$2H-\sqrt{4H^2-1}$. Actually, for $\tau=2H-\sqrt{4H^2-1}$, we find the
surface $\rho=\argth2H$ which is the vertical cylinder of constant mean
curvature $H$.

\begin{figure}[h]
\begin{center}
\resizebox{1\linewidth}{!}{\input{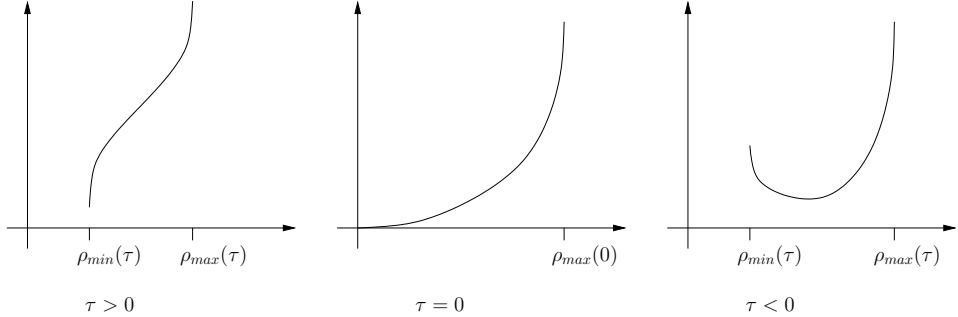}}
\caption{The solution of the equation}\label{fig:solution}
\end{center}
\end{figure}

For all values of $\tau<2H-\sqrt{4H^2-1}$, a solution $f$ can be defined for
$\rho\in [\rho_{min}(\tau),\rho_{max}(\tau)]$. This solution $f$ has the
following properties (see Figure~\ref{fig:solution}):
\begin{itemize}
\item if $0<\tau$, $0<\rho_{min}(\tau)$, $f$ is increasing and
$f'(\rho_{min}(\tau))=+\infty= f'(\rho_{max}(\tau))$.
\item if $\tau=0$, $0=\rho_{min}(0)$, $f$ is increasing, $f'(0)=0$ and $
f'(\rho_{max}(0))=+\infty$.
\item if $\tau<0$, $0<\rho_{min}(\tau)$, $f'(\rho_{min}(\tau))=-\infty$ and $
f'(\rho_{max}(\tau))=+\infty$.
\end{itemize}

In each case, the graph of the function $u$ is a piece of a cmc $H$ surface of
revolution that can be extended along its boundary by symmetry to produce a
complete rotationally invariant cmc $H$ surface $\boD_\tau$. When $\tau>0$ we
produce an embedded surface called unduloid. When $\tau=0$ we produce a cmc
sphere. For $\tau<0$, we get a non-embedded surface called a nodoid.

The parameter $\tau$ can be interpreted as a flux on the surface. More precisely,
if $\gamma$ is the circle $\boD_\tau\cap\{z=t\}$, $2\pi\tau$ is the flux of
$\boD_\tau$ along $\gamma$ in the direction $\partial_z$.

\section{The Alexandrov function and the main result}
\label{sec:alexandrov}
In this section we introduce the notion of Alexandrov functions. Then we study
the annular ends of a cylindrically bounded properly embedded cmc surface. We
then explain the proof of our main theorem.


\subsection{The Alexandrov function}\label{sec:alexfunc}
The notion of Alexandrov function was introduced by Korevaar, Kusner and Solomon
in \cite{KoKuSo} for $\R^3$. Let us explain what is the situation in
$\H^2\times\R$.

Let $\Gamma=(\gamma_t)_{t\in\R}$ be a smooth family of geodesic lines that
foliates $\H^2$. We define $\Gamma_t^+=\cup_{s>t}\gamma_s$, $\Gamma_t^+$ is a
half hyperbolic space bounded by $\gamma_t$; we also define
$\Gamma_t^-=\cup_{s<t}\gamma_s$. Let $\Pi_t$ be the vertical plane
$\gamma_t\times\R$ in $\H^2\times\R$. Let $S_t$ denote the symmetry of $\H^2$ with
respect to $\gamma_t$. $S_t$ extend to the symmetry of $\H^2\times \R$ with
respect to $\Pi_t$, we still denote by $S_t$ this symmetry. Among all these
foliations, we say that $\Gamma$ is a translation foliation if all the geodesic
lines $\gamma_s$ are orthogonal to one common geodesic line $g$.

For an open interval $I$ of $\R$, let $G$ be a cylindrically bounded domain
in $\H^2\times \bar I$ such that $\partial G\cap(\H^2\times I)$ is a smooth
connected surface $\Sigma$ with possibly non empty boundary in the horizontal
slices at height given by the end points of $I$.

Let us fix some $z_0$ in $I$, we focus on what happens at height
$z_0$. So we denote $\Sigma_{z_0}=\Sigma\cap \{z=z_0\}$ which bounds
$G_{z_0}=G\cap \{z=z_0\}$ and we consider that $\Gamma$ foliates
$\H^2\times\{z_0\}$. Since $G$ is cylindrically bounded, for $t$ large
$\Sigma_{z_0}(t)=\Sigma_{z_0}\cap \Gamma_t^+$ is empty. If $\Sigma_{z_0}$ is non
empty, there is a largest $t_1$ such that $\Sigma_{z_0}(t_1)\neq\emptyset$. For
$t\le t_1$, we consider $\widetilde{\Sigma}_{z_0}(t)=S_t(\Sigma_{z_0}(t))$ the
symmetric of $\Sigma_{z_0}(t)$ with respect $\gamma_t$. We notice that
$\widetilde{\Sigma}_{z_0}(t)$ is included in $\Gamma_t^-$ and that, for small
$t$, $\widetilde{\Sigma}_{z_0}(t)$ is outside $G_{z_0}$. So we can define 
$$
t_2(z_0)=\sup \{t\le t_1 \ |\,
\widetilde{\Sigma}_{z_0}(t)\cap\Sigma_{z_0}\neq\emptyset\}.
$$
We also define
$$
t_3(z_0)=\sup \{t\le t_1\ |\, \exists p\in \Sigma_{z_0}\cap \gamma_t,\
T_p\Sigma\text{ is orthogonal to }\gamma_t\times\R\}.
$$

Finally we define $\alpha_\Gamma(z_0)=\max(t_2(z_0),t_3(z_0))$, if
$\Sigma_{z_0}$ is empty we define $\alpha_\Gamma(z_0)=-\infty$. This number can
be understood as the first time where there is a contact between
$\widetilde{\Sigma}_{z_0}(t)$ and $\Sigma_{z_0}\cap \Gamma_t^-$. In fact, for
$t>\alpha_\Gamma(z_0)$, $\widetilde\Sigma_{z_0}(t)\subset G_{z_0}$. The function
$z_0\mapsto \alpha_\Gamma(z_0)$ is called the \emph{Alexandrov function} of
$\Sigma$ associated to the foliation $\Gamma$. This Alexandrov function has a
first important property.

\begin{lem}\label{lem:semicont}
The Alexandrov function $\alpha_\Gamma$ is upper semi-continuous.
\end{lem}
\begin{proof}
Let $(z_n)$ be a sequence in $I$ converging to $z_0$ in $I$. If $\Sigma_{z_0}$
is empty, $\Sigma_{z_n}$ is empty for large
$n$, so $\alpha_\Gamma$ is upper semi-continuous at $z_0$ if
$\alpha_\Gamma(z_0)=-\infty$. So we can assume $\alpha_\Gamma(z_n)>-\infty$ for
every $n$. We can also assume that either $\alpha_\Gamma(z_n)=t_2(z_n)$ for all
$n$
or $\alpha_\Gamma(z_n)=t_3(z_n)$ for all $n$. In the second case, we have a
sequence of point $p_n\in \Sigma_{z_n}\cap\gamma_{t_3(z_n)}$ such that
$T_{p_n}\Sigma$ is orthogonal to $\gamma_{t_3(z_n)}\times\R$. We can assume that
$t_3(z_n)\rightarrow \limsup \alpha_\Gamma(z_n)=t_0$ and $p_n\rightarrow p_0$.
Then $p_0\in\Sigma_{z_0}\cap \gamma_{t_0}$ and $T_{p_0}\Sigma$ is orthogonal to
$\gamma_{t_0}\times\R$. This implies that $t_0\le t_3(z_0)\le
\alpha_\Gamma(z_0)$.

If we are in the first case, there is a sequence of point $p_n$ in
$\Sigma_{z_n}(t_2(z_n))$ such that $S_{t_2(z_n)}(p_n)\in \Sigma_{z_n}$. We can
assume that $t_2(z_n)\rightarrow \limsup \alpha_\Gamma(z_n)=t_0$ and
$p_n\rightarrow p_0$. Thus $p_0\in\Sigma_{z_0}(t_0)$ or $p_0\in\Sigma_{z_0}\cap
\gamma_{t_0}$. In the first case $S_{t_0}(p_0)\in\Sigma_{z_0}$ so $t_0\le
t_2(z_0)\le \alpha_\Gamma(z_0)$. In the second, since $S_{t_2(z_n)}(p_n)\in
\Sigma_{z_n}$ and converge to $p_0$, the tangent space $T_{p_0}\Sigma$ is
orthogonal to $\gamma_{t_0}\times\R$. So $t_0\le t_3(z_0)\le
\alpha_\Gamma(z_0)$. This finishes the proof.
\end{proof}

The second important property is a consequence of the maximum principle when
$\Sigma$ has constant mean curvature.

\begin{lem}\label{lem:nomax}
Assume that $\Sigma$ is connected, has constant mean curvature and the mean
curvature vector points into $G$. If the Alexandrov function $\alpha_\Gamma$ has
a local maximum at $z$, $\Sigma$ is symmetric with
respect to $\Pi_{\alpha_\Gamma(z)}$.
\end{lem}
\begin{proof}
Assume that $z$ is a local maximum of $\alpha_\Gamma$. Let $p\in\Sigma_z$ such
that $p\in\Sigma_z(\alpha_\Gamma(z))$ and $S_{\alpha_\Gamma(z)}(p)\in\Sigma_z$
if $\alpha_\Gamma(z)=t_2(z)$ or $p\in\gamma_{\alpha_\Gamma(z)}\times\R$ and
$T_p\Sigma$ is orthogonal to $\gamma_{\alpha_\Gamma(z)}\times\R$ if
$\alpha_\Gamma(z)=t_3(z)$. In both case, since $z$ is a local maximum of
$\alpha_\Gamma$, near $S_{\alpha_\Gamma(z)}(p)$, $S_{\alpha_\Gamma(z)}(\Sigma)$
is on one side of $\Sigma$. Moreover these two surfaces have the same
orientation at $S_{\alpha_\Gamma(z)}(p)$. So applying the maximum principle or
the maximum principle at the boundary for cmc surface, we get that $\Sigma$ is
symmetric with respect to $\Pi_{\alpha_\Gamma(z)}$
\end{proof}

These two results has the following consequence

\begin{lem}\label{lem:monot}
Let $\Sigma$ be as in Lemma~\ref{lem:nomax}, if $\alpha_\Gamma>-\infty$
on $[a,b]$ then:
\begin{itemize}
\item $\alpha_\Gamma$ is monotonous, or
\item there exists $c\in [a,b]$ such that $\alpha_\Gamma$ is decreasing on $[a,c)$ and increasing on $(c,b]$.
\end{itemize}
\end{lem}
\begin{proof}
First, from Lemma~\ref{lem:nomax}, we remark that, if $\alpha_\Gamma$ has a
local maximum at $z$, $\Sigma$ is symmetric with respect to
$\Pi_{\alpha_{\Gamma}(z)}$ so $\alpha_\Gamma(z)$ is a global minimum for
$\alpha_\Gamma$ on $[a,b]$ so $\alpha_\Gamma$ is constant close to $z$.

Let $x,y$ be in $[a,b]$. Since $\alpha_\Gamma$ is upper semi-continuous, $\sup_{[x,y]}\alpha_\Gamma$ is reached somewhere in $[x,y]$. If it is in the inside, the above remark implies that $\alpha_\Gamma$ is constant on $[x,y]$. Thus the maximum is always reached
at $x$ or $y$. So we have
$\sup_{[x,y]}\alpha_\Gamma=\sup\{\alpha_\Gamma(x),\alpha_\Gamma(y)\}$.

Let $(x_n)$ be a monotonous sequence converging to $c\in[a,b]$ such that $\lim
\alpha_\Gamma (x_n)=\inf_{[a,b]}\alpha_\Gamma$. We assume that $(x_n)$ decreases
(the same argument can be done if it increases). 

Let us consider
$ x< y< c$. If $\alpha_\Gamma( y)=\inf_{[a,b]}\alpha_\Gamma$, for any $z \in
(y,c)$, $\alpha_\Gamma(z)\le \alpha_\Gamma(x_n)$ so
$\alpha_\Gamma(z)=\inf_{[a,b]}\alpha_\Gamma$. $\alpha_\Gamma$ is constant on
$(y,c)$ so it is constant on $[a,b]$. So we can assume
$\alpha_\Gamma(y)>\lim\alpha_\Gamma(x_n)$. Since $\alpha_\Gamma(y)\le
\max\{\alpha_\Gamma(x),\alpha_\Gamma(x_n)\}$, we get $\alpha_\Gamma(x)\ge
\alpha_\Gamma(y)$: $\alpha_\Gamma$ decreases on $[a,c)$.

Let us now consider $c<x<y$. As above we can assume
$\alpha_\Gamma(x)>\lim\alpha_\Gamma (x_n)$. Since $\alpha_\Gamma(x)\le
\max\{\alpha_\Gamma(y),\alpha_\Gamma(x_n)\}$, we get $\alpha_\Gamma(x)\le
\alpha_\Gamma(y)$: $\alpha_\Gamma$ increases on $(c,b]$.

When $c=a$ or $b$, $\alpha_\Gamma$ can be monotonous.
\end{proof}

If $\Sigma$ is a weakly embedded cmc surface which is cylindrically bounded, we
notice that the Alexandrov functions can also be defined on it. These functions
satisfies also to Lemmas~\ref{lem:nomax} and \ref{lem:monot}.


\subsection{Asymptotical Delaunay ends} 

In this subsection we prove that a cylindrically bounded annular ends of a cmc
surface is asymptotic to a Delaunay surface $\boD_\tau$.

Let $A$ be a properly embedded annular end with cmc $H$ and which is
cylindrically bounded. The annulus $A$ can be viewed as the punctured disk
embedded in $\H^2\times \R$ with boundary in $\H^2\times\{0\}$. The height
function converges to $\pm\infty$ at the puncture. If the limit of the height
function is $+\infty$ $A$ is said to be a top end and if the limit is $-\infty$
$A$ is a bottom end. In the sequel we always study top ends, since the
properties of bottom ends can be deduced by symmetry with respect to $\{z=0\}$.
The annular end $A$ bounds a cylindrically bounded domain $G$ so we can define
Alexandrov functions on $A$.

From Proposition~\ref{prop:linarea}, $A$ has linear area growth.
Proposition~\ref{prop:boundcurv} gives a uniform bound on the second
fundamental form on $A$. This two properties are sufficient to study the limit
of slight-back sequences $t_{-z_n}(A)$ where $(z_n)$ is an increasing sequence
going to $+\infty$ and $t_{\bar z}$ denote the vertical translation by $\bar z$.

\begin{prop}\label{prop:asymp}
Let $A$ be a properly embedded annular end with cmc $H$ which is cylindrically
bounded ($H>1/2$). Let $(z_n)$ be an increasing sequence converging to $+\infty$.
There is a parameter $\tau\in (0,2H-\sqrt{4H^2-1}]$ that depends only on $A$
(not on $(z_n)$) and a
subsequence $(z_{n'})$ of $(z_n)$ such that $t_{-z_{n'}}(A)$ converges to a
rotationally invariant Delaunay surface $\boD_\tau$. Moreover the axis of
$\boD_\tau$ only depends on $A$.
\end{prop}

\begin{proof}
Let $s\mapsto\gamma(s)$ be a geodesic in $\H^2$. Let $\gamma_s^1$ denote the
geodesic line of $\H^2$ orthogonal to $\gamma$ at $\gamma(s)$. Then
$\Gamma^1=(\gamma_s^1)_{s\in\R}$ is a foliation of $\H^2$, so we can consider
the Alexandrov function $\alpha_{\Gamma^1}$ of $A$. From Lemma~\ref{lem:monot},
$\alpha_{\Gamma^1}$ is monotonous close to $+\infty$. Moreover, it is bounded so
it has a limit at $+\infty$. by changing the parametrization of $\gamma$, we
assume that this limit is $0$.

Let $\theta$ be an irrational angle and consider the geodesic
line $\gamma_t^2$ of $\H^2$ which meets $\gamma_0^1$ with an angle $\theta$ at
$\gamma_0^1(t)$. $\Gamma^2=(\gamma_t^2)_{t\in\R}$ is a foliation of $\H^2$ so we
consider the Alexandrov function $\alpha_{\Gamma^2}$ of $A$. As above, this
function has a limit at $+\infty$ and we can assume it is $0$. We denote by $p$
the point where $\gamma_0^1$ and $\gamma_0^2$ meet. We notice that the position
of this point $p$ will fix the axis of the Delaunay limit surface.

Now let us consider an increasing sequence $(z_n)$ with limit $+\infty$. From
Subsection~\ref{sec:conv}, a subsequence of $t_{-z_n}(A)$ (still denoted
$t_{-z_n}(A)$) converges to a properly weakly embedded surface $\Sigma$ with
constant mean curvature $H$. Moreover $\Sigma$ is cylindrically bounded. 

\begin{claim}
The surface $\Sigma$ is connected and non-compact.
\end{claim}
\begin{proof}[Proof of the Claim]
Assume $\Sigma$ has a compact connected component $\Sigma'$. This implies that,
for large $n$, there is a part of $t_{-z_n}(A)$ that is graph over $\Sigma'$. So
$A$ would possess a compact component, this gives a contradiction.

Assume now that $\Sigma$ has two non compact connected components $\Sigma'$ and
$\Sigma''$. We recall that the height function on $\Sigma'$ and $\Sigma''$ can
not be lower or upper bounded (Lemma~\ref{lem:heightbound}). So there is a
connected component of $\Sigma'\cap\{0\le z\le M\}$ and one of
$\Sigma''\cap\{0\le z\le M\}$ with boundary in both $\{z=0\}$ and $\{z=M\}$. So
this implies that for large $n$,
$A\cap\{z_n\le z\le z_n+M\}$ possesses at least two connected components with
boundary in both $\{z=z_n\}$ and $\{z=z_n+M\}$. This is in contradiction with
Lemma~\ref{lem:compannulus} when $M$ is large.
\end{proof}

The existence of a limit for the Alexandrov function $\alpha_{\Gamma^1}$ and $\alpha_{\Gamma^2}$ has the following consequence.

\begin{claim}
The surface $\Sigma$ is symmetric with respect to $\gamma_0^1\times\R$ and
$\gamma_0^2\times\R$.
\end{claim}

\begin{proof}[Proof of the Claim]
We only write the proof for $\gamma_0^1\times\R$. Let $S_t^1$ denote the
symmetries associated to the foliation $\Gamma^1$. We want to prove that
$S_0^1(\Sigma\cap\Gamma_0^+)$ is on one side of $\Sigma\cap\Gamma_0^-$ and these
surfaces touch each other (eventually in the boundary). If it is true the
maximum principle applies and we get the symmetry with respect to
$\gamma_0^1\times\R$.

So let us fix $z_0$ in $\R$ a regular value of the height function on $\Sigma$.
Because of the properness, any value $t$ close to $z_0$ is also a regular value
of the height function. So $\Sigma\cap \{z_0-\eps\le z \le z_0+\eps\}$ consists
in a finite union of annuli transverse to horizontal slices of $\H^2\times\R$.

Let $A_n$ denote the sequence $t_{-z_n}(A)$ and $G_n=t_{-z_n}(G)$ (where $G$ is
the domain bounded by $A$). We also denote by $\overline G$ the domain bounded
by $\Sigma$ and $t_n=\alpha_{\Gamma^1}(t+z_n)$ for some
$t\in(z_à-\eps,z_0+\eps)$: we have $t_n\rightarrow 0$. Because of the
transversality and the convergence, we have $A_n\cap \{z=t\}\rightarrow
\Sigma\cap \{z=t\}$.

Besides, from the definition of the Alexandrov function, $S_{t_n}^1(A_n\cap
\Gamma_{t_n}^+\cap\{z=t\})$ is in $G_n\cap\{z=t\}$. Thus, by taking the limit,
$S_0^1(\Sigma\cap \Gamma_0^+\cap\{z=t\})$ is in $\overline G\cap\{z=t\}$. Since
we have this for any $t\in(z_0-\eps,z_0+\eps)$, we have
$S_0^-(\Sigma\cap \Gamma_0^+\cap \{z_0-\eps<z<z_0+\eps\})\subset G$ so it is on
one side of $\Sigma\cap \Gamma_0^-$. Besides
as in the proof of Lemma~\ref{lem:semicont}, there is a contact point between
these two surfaces at every height $t\in(-z_0-\eps,z_0+\eps)$. Thus the maximum
principle applies to prove that $\Sigma$ is invariant by $S_0^1$ (here we use
the connectedness of $\Sigma$).
\end{proof}




Since $\Sigma$ is symmetric with respect to $\gamma_0^1\times\R$ and $\gamma_0^2\times\R$, it is invariant by the
rotation $R$ of angle $2\theta$ around the vertical axis $\{p\}\times\R$. Since
$2\theta$ is irrational, $\Sigma$ is invariant by rotation around
$\{p\}\times\R$.

This implies that $\Sigma$ is equal to a Delaunay surface $\boD_\tau$ of axis
$\{p\}\times\R$. The
height $0$ is a regular value of $\boD_\tau$ and $2\pi\tau$ is the flux of
$\boD_\tau$ along $\boD_\tau\cap \{z=0\}$ in the direction $\partial_z$. So
$2\pi\tau$ is the limit of the flux of $A$ along $A\cap\{z=z_n\}$ in the
direction $\partial_z$. But this flux does not depends on $n$ and is equal to the
flux of $A$ along $A\cap \{z=0\}$ in the direction $\partial_z$. So the
parameter $\tau$ only depends on $A$.
\end{proof}


\subsection{The main theorem}

In this section, we settle the main theorem of this paper and explain its proof.

Actually, the main theorem is based on the following proposition that will be proved
in Section~\ref{sec:monotonous}.

\begin{prop}\label{prop:decrease}
Let $A$ be a properly embedded annular top end with cmc $H$ which is
cylindrically bounded. Let $\Gamma=(\gamma_s)_{s\in\R}$ be a translation foliation of
$\H^2$ by geodesic lines. The Alexandrov function $\alpha_\Gamma
:\R^+\rightarrow \R$ is then decreasing.
\end{prop}

With this proposition we can prove our main result.

\begin{thm}\label{th:main}
Let $\Sigma$ be a properly embedded cmc surface in $\H^2\times\R$. If $\Sigma$
has finite topology and is cylindrically bounded, $\Sigma$ is a Delaunay surface
(\textit{i.e.} $\Sigma$ is rotationally invariant).
\end{thm}

If $\Sigma$ is compact the result is already known \cite{HsHs}, so we focus on
the non compact case.
\begin{proof}
Let us consider $\Gamma=(\gamma_s)_{s\in\R}$ be a translation foliation of
$\H^2$. Let us denote by $E_1^+,\dots,E_p^+$ the annular top ends of $\Sigma$ and
$E_1^-,\dots,E_q^-$ the annular bottom ends of $\Sigma$. We consider the
Alexandrov functions $\alpha_{\Gamma,\Sigma}$ and $\alpha_{\Gamma,E_i^\pm}$. We
can assume that the Alexandrov functions $\alpha_{\Gamma,E_i^+}$ are defined on
$[M,+\infty)$ and the $\alpha_{\Gamma,E_i^-}$ are defined on $(-\infty,-M]$. By
Proposition~\ref{prop:decrease}, the functions $\alpha_{\Gamma,E_i^+}$ decrease
and the functions $\alpha_{\Gamma,E_i^-}$ increase.

Besides, on $[M,+\infty)$, we have:
$$
\alpha_{\Gamma,\Sigma}(z)=\max_{1\le i\le p}\alpha_{\Gamma,E_i^+}(z)
$$
and, on $(-\infty,-M]$, we have:
$$
\alpha_{\Gamma,\Sigma}(z)=\max_{1\le i\le q}\alpha_{\Gamma,E_i^-}(z)
$$
So the function $\alpha_{\Gamma,\Sigma}$ increases on $(-\infty,-M]$ and
decreases on $[M,+\infty)$. By Lemma~\ref{lem:monot}, this implies that
$\alpha_{\Gamma,\Sigma}$ is constant and $\Sigma$ is symmetric with respect
to some $\gamma_s\times\R$.

Let $\Gamma^1$ be a translation foliation of $\H^2$, $\Sigma$ is then symmetric
with respect to some $\gamma_s^1$. We can assume that it is symmetric with
respect to $\gamma_0^1\times \R$. Let $\Gamma_2$ be the translation foliation of
$\H^2$ composed by the geodesic line orthogonal to $\gamma_0^1$. $\Sigma$ is
then symmetric with respect to some $\gamma_s^2\times\R$. We can also assume it
is $\gamma_0^2$. Let $p$ be the intersection point of $\gamma_0^1$ and
$\gamma_0^2$. Let $g$ be a geodesic passing by $p$ and $\Gamma^3$ be the
translation foliation composed by the geodesic lines orthogonal to $g$. $\Sigma$
is then symmetric with respect to some $\gamma_s^3\times\R$. Since $\Sigma$ is
cylindrically bounded $\gamma_s^3$ passes by $p$. This implies that $\Sigma$ is
symmetric with respect to any vertical plane passing by $p$, so $\Sigma$ is
invariant by rotation around the vertical axis $p\times\R$.
\end{proof}

\section{Horizontal Killing graphs}
\label{sec:horigraph}
Let us consider a new model for $\H^2$: $\H^2=\{(s,r)\in\R^2\}$ with the metric
$\dd r^2+ (\cosh r)^2\dd s^2$. In this model $\{r=0\}$ is a geodesic,
$\{r=c\}$ are its equidistant lines and $\{s=c\}$ are the geodesic lines
orthogonal to $\{r=0\}$. Moreover $\partial_s$ is the Killing
vector field corresponding to the translation along $\{r=0\}$. In this section,
we will use this model of $\H^2$ to describe $\H^2\times\R$. The surfaces
$\{s=c\}$ are then totally geodesic flat planes.

Let $\Ome$ be a domain in $\R^2$ and $u$ be a smooth function on $\Ome$. Using
the above model for $\H^2$, we can consider the surface in $\H^2\times\R$
parametrized by $(r,z)\mapsto (u(r,z),r,z)$. Such a surface is called the
\emph{horizontal Killing graph} of $u$, it is transverse to the Killing vector
field $\partial_s$ and any integral curve of $\partial_s$ intersect at most once
the surface.


\subsection{The mean curvature equation}
In the following, we are interested in horizontal Killing graphs with constant
mean curvature $H_0$. This condition implies that the function $u$ satisfies a
partial differential equation.

Let $\Ome$ and $u$ be as above. A unit normal vector to the horizontal Killing
graph of $u$ is given by the following expression:
\begin{equation}\label{eq:normal}
N=\frac{-\partial_s+(\cosh r)^2\nabla u}{\cosh r\sqrt{1+ (\cosh r)^2 |\nabla
u|^2}}
\end{equation}
where $\nabla$ is the Euclidean gradient operator and $|\cdot|$ is the Euclidean
norm. In the sequel, we will use this unit normal vector to compute the mean
curvature of a horizontal Killing graph.

\begin{lem}
Let $\Ome$ and $u$ be as above, the mean curvature $H$ of the horizontal Killing
graph of $u$ satisfies:
\begin{equation}
-2H\cosh r=\Div\frac{(\cosh r)^2\nabla u}{\sqrt{1+(\cosh r)^2|\nabla u|^2}}
\end{equation}
with $\Div$ the Euclidean divergence operator.
\end{lem}

\begin{proof}
We extend the vector field $N$ to the whole $\R\times\Ome$ by using the
expression given in \eqref{eq:normal}. The mean curvature of the horizontal
Killing graph of $u$ is then given by 
\begin{align*}
-2H&=\Div_{\H^2\times\R}N\\
&=(\bnabla_{\frac{\partial_s}{\cosh r}}N, \frac{\partial_s}{\cosh r})+
(\bnabla_{\partial_r}N, \partial_r)+ (\bnabla_{\partial_z}N, \partial_z)
\end{align*}
Let $W$ denote $\sqrt{1+(\cosh r)^2|\nabla u|^2}$, we then have:
\begin{align*}
(\bnabla_{\frac{\partial s}{\cosh r}}N, \frac{\partial_s}{\cosh r})&=
\frac{1}{(\cosh r)^2}\left(\frac{-1}{W\cosh r }(\bnabla_{\partial_s} \partial_s,
\partial_s)+ \frac{\cosh
r}{W}(\bnabla_{\partial_s}\nabla u,\partial_s)\right)\\
&=\frac{1}{W\cosh r}(\bnabla_{\nabla u}\partial_s,\partial_s)\\
&=\frac{\nabla u\cdot\cosh r}{W}
\end{align*}
and for $a=r$ or $a=z$:
\begin{align*}
(\bnabla_{\partial_a}N, \partial_a)&=(\bnabla_{\partial_a}\frac{-\partial_s}{W\cosh r}, \partial_a )+ (\nabla_{\partial_a}\frac{(\cosh r)\nabla u}{W},\partial_a)\\
&=(\nabla_{\partial_a}\frac{(\cosh r)\nabla u}{W},\partial_a)
\end{align*}
Summing all these terms, we get:
\begin{align*}
-2H&=\frac{\nabla u\cdot\cosh r}{W}+\Div(\frac{(\cosh r)\nabla u}{W})\\
&=\frac{1}{\cosh r}\left(\frac{(\nabla(\cosh r),(\cosh r)\nabla u)}{W}+ (\cosh
r)\Div(\frac{(\cosh r)\nabla u}{W})\right)\\
&=\frac{1}{\cosh r}\Div(\frac{(\cosh r)^2\nabla u}{W})
\end{align*}
\end{proof}

Thus the constant mean curvature $H_0$ equation for a function $u$ can be written
\begin{equation}\label{eq:cmc1}
\Div\frac{(\cosh r)^2\nabla u}{\sqrt{1+(\cosh r)^2|\nabla u|^2}}=-2H_0\cosh r
\end{equation}
or after expanding all the terms
{\footnotesize
\begin{multline}\label{eq:cmc2}
\left(\frac{1+(\cosh r)^2 u_z^2}{2+(\cosh r)^2|\nabla u|^2}\right)u_{rr}-
2\left(\frac{(\cosh r)^2 u_ru_z}{2+(\cosh r)^2|\nabla u|^2}\right)u_{rz}
+\left(\frac{1+(\cosh r)^2 u_r^2}{2+(\cosh r)^2|\nabla
u|^2}\right)u_{zz}\\+(\tanh r) u_r=-2H_0 \frac{(1+(\cosh r)^2|\nabla
u|^2)^{3/2}}{(\cosh r)(2+(\cosh r)^2|\nabla u|^2)}
\end{multline}}

For $H_0=0$, we get the minimal surface equation:
\begin{equation}\label{eq:min}
\Div\frac{(\cosh r)^2\nabla u}{\sqrt{1+(\cosh r)^2|\nabla u|^2}}=0
\end{equation}

We notice that the maximum principle is true for these equations. Thus we have
uniqueness of a solution to the Dirichlet problem associated to these equations
on bounded domains.


\subsection{A gradient estimate}

An important result concerning solutions of \eqref{eq:cmc1} is the following
gradient estimate.

\begin{prop}\label{prop:estim}
Let $u$ be a nonnegative solution of \eqref{eq:cmc1} on a disk centered at
$p=(r_p,z_p)$ and radius $R$. Then there is a constant $M$ that depends only on
$r_p$, $R$ and $H_0$ such that 
$$
|\nabla u|(p)\le \max(2,32M(u(p)/R))e^{6Mu(p)+4M^2(u(p)/R)^2}
$$
\end{prop}

The proof of this result is similar to the one of the gradient estimate proved
by J.~Spruck in \cite{Spr}; but our result does not seem to be a corollary of
his result.

Before beginning the proof, let us make some preliminary computation. So let $u$
be as in the proposition and let $\Sigma$ denote the horizontal Killing
graph of $u$. We denote $\tilde N=-N$ (see \eqref{eq:normal}) and define
$\nu=(\tilde N,\partial_s)$ and $\mu=\nu/(\cosh r)$. We have $\nu>0$ and, $\partial_s$ being a Killing vector field,
$$
\Delta_\Sigma\nu=-(Ric(\tilde N,\tilde N)+|A|^2)\nu
$$
where $|A|^2$ is the square of the norm of the second fundamental form and $Ric$
is the Ricci tensor.

Let us denote by $h$ the restriction of $s$ along $\Sigma$. We have
$$
\nabla_\Sigma h=\frac{1}{(\cosh r)^2}\partial_s^\top\quad \text{and}\quad
|\nabla_\Sigma h|^2=\frac{1}{(\cosh r)^2}(1-\mu^2)
$$
If $(e_1,e_2)$ is an orthonormal basis of $T\Sigma$ we have
\begin{align*}
\Delta_\Sigma h&=\Div_\Sigma(\frac{1}{(\cosh r)^2}\partial_s^\top)=-\frac{2\tanh
r}{(\cosh r)^2}(\partial_r^\top,\partial_s^\top)+\frac{1}{(\cosh
r)^2}\Div_\Sigma(\partial_s-\nu \tilde N)\\
&=\frac{2\tanh r}{\cosh r}\mu(\partial_r,\tilde N)+\frac{1}{(\cosh
r)^2}(-\nu)\sum_{i=1}^2(\bnabla_{e_i}\tilde N,e_i)\\
&=\frac{2\tanh r}{\cosh r}\mu(\partial_r,\tilde N)-\frac{2H_0\mu}{\cosh r}
\end{align*}

Let us define the distance function $d=((r-r_p)^2+(z-z_p)^2)^{1/2}$. The vector
field $\partial_d=((r-r_p)\partial_r+(z-z_p)\partial_z)/d$ is well defined in
$\H^2\times\R$ outside $\R\times \{p\}$ and has unit length; $d\partial_d$ is
well defined everywhere. We have:
$$
\nabla_\Sigma d^2=2d\partial_d^\top\quad \text{and}\quad |\nabla_\Sigma
d^2|^2=4d^2|\partial_d^\top|^2
$$
We denote $\tilde r=r-r_p$ and $\tilde z=z-z_p$. We then have:
\begin{align*}
\Delta_\Sigma d^2&=2\Div_\Sigma(\tilde r\partial_r^\top + \tilde
z\partial_z^\top)\\
&=2(|\partial_r^\top|^2+|\partial_z^\top|^2)+2\sum_{i=1}^2\Big(\tilde
r(\bnabla_{e_i}(\partial_r-(\partial_r,\tilde N)\tilde N,e_i)+ \tilde
z(\bnabla_{e_i}(\partial_z-(\partial_z,\tilde N)\tilde N,e_i)\Big)\\
&=2(1+\mu^2)+2\sum_{i=1}^2\Big(\tilde r(\bnabla_{e_i}\partial_r,e_i)-\big(\tilde
r(\partial_r,\tilde N)+\tilde z(\partial_z,\tilde N)\big)(\bnabla_{e_i}\tilde
N,e_i)\Big)\\
&=2(1+\mu^2)+2\tilde
r\sum_{i=1}^2(\bnabla_{e_i}\partial_r,e_i)-2H_0(d\partial_d,\tilde N)
\end{align*}
We define $f_1=\partial_s/(\cosh r)$, $f_2=\partial_r$ and $f_3=\partial_z$ an
orthonormal basis of $T\H^2\times\R$ and we write $e_i=\sum_j\lambda_i^jf_j$. We
then have
\begin{align*}
\sum_{i=1}^2(\nabla_{e_i}\partial_r,e_i)&= \sum_{i=1}^2\sum_{k,l=1}^3
\lambda_i^k\lambda_i^l(\bnabla_{f_k}\partial_r,f_l)\\
&=\sum_{i=1}^2(\lambda_i^1)^2\frac{1}{(\cosh r)^2}(\bnabla_{\partial_s}
\partial_r,\partial_s)\\
&=\tanh r(1-\mu^2)
\end{align*}
So $\Delta_\Sigma d^2=2(1+\mu^2)+2\tilde r\tanh r(1-\mu^2)-2H_0d(\partial_d,
\tilde N)$.
Using the above computations, we are ready to write the proof.

\begin{proof}[Proof of Proposition~\ref{prop:estim}]
Let us introduce the second order operator $Lf=\Delta_\Sigma f-
2\nu(\nabla_\Sigma\frac{1}{\nu}, \nabla_\Sigma f)$ on $\Sigma$. We notice that
the maximum principle is true for $L$. We have:
$$
\Delta_\Sigma\frac{1}{\nu}=\Div_\Sigma(-\frac{1}{\nu^2}\nabla_\Sigma\nu)=
2\nu|\nabla_\Sigma\frac{1}{\nu}|^2+(Ric(\tilde N,\tilde N)+|A|^2)\frac{1}{\nu}
$$
Since $Ric(\tilde N,\tilde N)\ge -1$, we have $L\frac{1}{\nu}\ge-\frac{1}{\nu}$.
Let us define $v=\eta\frac{1}{\nu}$ with $\eta$ a positive function. We have:
$$
Lv=\eta L\frac{1}{\nu}+\frac{1}{\nu}\Delta_\Sigma\eta\ge
(\Delta_\Sigma\eta-\eta)\frac{1}{\nu}
$$
We define on $\Sigma$ the function
$\phi=((-\frac{h}{2h_0}+1-\eps-(\frac{d}{R})^2)^+$ which
is less than $1$ ($\eps>0$) where $h_0=u(p)=h(P)$ with $P=(u(p),p)$. Moreover,
$\phi(P)=1/2-\eps$ and $\phi=0$ close to $\partial\Sigma$. We define
$\eta=e^{K\phi}-1$ with $K$ a positive constant that will be chosen later. We
then have $\max v>0$ and it is reached inside the support of $\phi$.

We have $\Delta_\Sigma\eta= \Div_\Sigma(e^{K\phi}K\nabla_\Sigma\phi)=
e^{K\phi}(K^2|\nabla_\Sigma\phi|^2+K\Delta_\Sigma\phi)$ so:
\begin{align*}
\Delta_\Sigma\eta-\eta&=e^{K\phi}(K^2|\nabla_\Sigma\phi|^2+K\Delta_\Sigma\phi-1)
+1\\
&\ge e^{K\phi}(K^2|\nabla_\Sigma\phi|^2+K\Delta_\Sigma\phi-1)
\end{align*}
We have
\begin{align*}
|\nabla_\Sigma \phi|^2&=|-\frac{1}{2h_0}\nabla_\Sigma h-\frac{1}{R^2}
\nabla_\Sigma d^2|^2\\
&=\frac{1}{4h_0^2\cosh^2 r}(1-\mu^2)+\frac{4d^2}{R^2}|\partial_d^\top|^2+
\frac{2d}{h_0R^2\cosh^2r}(\partial_s^\top,\partial_d^\top)\\
&\ge \frac{1}{4h_0^2\cosh^2 r}(1-\mu^2)-\frac{2d}{h_0R^2\cosh r}\mu(\partial_d,N)\\
&\ge\frac{1}{4h_0^2\cosh^2 r}(1-\mu^2-8\frac{h_0}{R}M\mu)
\end{align*}
where we use $d\le R$ and $M$ is a constant chosen to be larger than $\cosh
r\sqrt{4+2R+2H_0R}$ on the disk of center $p$ and radius $R$. So:
$$
\text{ if }\mu\le \min(\frac{1}{2},\frac{R}{32Mh_0}),\quad
|\nabla_\Sigma\phi|^2\ge\frac{1}{8h_0^2(\cosh r)^2}
$$
Besides, we have
\begin{align*}
\Delta_\Sigma\phi&=-\frac{1}{2h_0}\Delta_\Sigma h-\frac{1}{R^2}\Delta_\Sigma
d^2\\
&=-\frac{1}{2h_0}(\frac{2\tanh r}{\cosh r}\mu(\partial_r,\tilde N)-
\frac{2H_0\mu}{\cosh r})-\frac{1}{R^2}(2(1+\mu^2)+2\tilde r\tanh
r(1-\mu^2)-2H_0d(\partial_d,\tilde N))\\
&\ge-\frac{\mu}{h_0\cosh r}-\frac{1}{R^2}(2(1+\mu^2)+2R(1-\mu^2)+2H_0R)\\
&\ge-\frac{1}{h_0^2(\cosh r)^2}\Big(Mh_0+\frac{h_0^2(\cosh
r)^2}{R^2}(4+2R+2H_0R)\Big)\\
&\ge-\frac{1}{h_0^2(\cosh r)^2}\Big(Mh_0+\frac{h_0^2M^2}{R^2}\Big)
\end{align*}
We deduce from the above computation that, if $\mu\le
\min(\frac{1}{2},\frac{R}{32Mh_0})$,
\begin{align*}
K^2|\nabla_\Sigma\phi|^2+K\Delta_\Sigma\phi-1 &\ge \frac{1}{8h_0^2(\cosh
r)^2}K^2 +K(-\frac{1}{h_0^2(\cosh r)^2})\Big(Mh_0+\frac{h_0^2M^2}{R^2}\Big)-1\\
&\ge \frac{1}{8h_0^2(\cosh
r)^2}\Big(K^2-8K(Mh_0+\frac{h_0^2M^2}{R^2})-8h_0^2M^2\Big)
\end{align*}
So if $K=(12Mh_0+8h_0^2M^2/R^2)$ we obtain that $K^2|\nabla_\Sigma\phi|^2+
K\Delta_\Sigma\phi-1>0$ and then $Lv>0$. By the maximum principle applied to
$L$, it implies that the maximum of $v$ can only be attained at a point $q$
where $\mu\ge\min(\frac{1}{2},\frac{R}{32Mh_0})$. This implies that
$$
v(p)=(e^{K(1/2-\eps)}-1)\frac{1}{\nu(p)}\le \frac{e^K-1}{\nu(q)}\le
\frac{e^K-1}{\min(\frac{1}{2},\frac{R}{32Mh_0})}
$$
So letting $\eps$ tending to $0$ we get:
$$
\nu(p)\ge \min(\frac{1}{2},\frac{R}{32Mh_0})e^{-K/2}
$$
So:
$$
|\nabla u|(p)\le \max(2,32M(h_0/R))e^{6Mh_0+4M^2(h_0/R)^2}
$$
\end{proof}


\subsection{An existence result for the Dirichlet problem}

In this subsection, we give a result about the existence of a solution of
the Dirichlet problem for the equation \eqref{eq:cmc1} on small domains.
Actually, it is a consequence of the work of J.~Serrin in \cite{Ser}.

\begin{prop}\label{prop:exist}
Let $p=(r_p,z_p)$ be a point $\R^2$ and $H_0$ be a nonnegative constant. Then,
there exists a constant $R>0$ that depends only on $H_0$ and $|r_p|$ such that
the Dirichlet problem for the equation \eqref{eq:cmc1} can be solved on the disk
$D(p,\tilde R)$ centered at $p$ and radius $\tilde R$ less than $R$. More
precisely, for any continuous function $\phi$ on the boundary of $D(p,\tilde R)$
($\tilde R\le R$) there exists $u\in C^2(D(p,\tilde R))\cap
C^0(\overline{D(p,\tilde R)})$ such that $u$ solves \eqref{eq:cmc1} and $u=\phi$
on the boundary of the disk.
\end{prop}

\begin{proof}
if $\phi$ is $C^2$, the result is a consequence of Theorem~14.3 in \cite{Ser}. We
notice that the hypotheses of this theorem are satisfied by
Equation~\eqref{eq:cmc1}. In fact in order to have the same notation as
J.~Serrin, the equation has to be written in the form \eqref{eq:cmc2}. Moreover,
since the coefficients of \eqref{eq:cmc2} only depend on $r$ and $H_0$, the
radius $R$ only depends on $r_p$ and $H_0$.

When $\phi$ is only continuous, we proceed by approximation. Let $(\phi_n)$ be a
sequence of $C^2$ functions converging to $\phi$ in the $C^0$ norm. We denote by
$u_n$ the solution of \eqref{eq:cmc1} with $\phi_n$ as boundary value. The
sequence $u_n$ is uniformly bounded. So, by the gradient estimate
(Proposition~\ref{prop:estim}) and elliptic estimates, the sequence $(u_n)$
converges to a solution $u$
of \eqref{eq:cmc1} on the disk. Let us consider $\eps>0$ and $n\in\N$ such that
$\phi_n-\eps\le \phi-\eps/2<\phi+\eps/2\le \phi_n+\eps$, for $m$ large we have
$\phi_n-\eps\le \phi_m\le \phi_n+\eps$. So by the maximum principle,
$u_n-\eps\le u_m\le u_n+\eps$. This implies $u_n-\eps\le u\le u_n+\eps$, so on
the boundary of the disk
$$
\phi-3\eps/2\le \phi_n-\eps\le \liminf_{\partial D(p,\tilde R)}u\le \limsup_{\partial
D(p,\tilde R)} u\le \phi_n+\eps\le \phi+3\eps/2
$$
Letting $\eps$ going to $0$, we see that $u$ is continuous up to the boundary
and $u=\phi$ there.
\end{proof}


\subsection{A uniqueness result}
In this section we give a uniqueness result for the Dirichlet problem associated
to \eqref{eq:cmc1} when the domain is unbounded.

\begin{prop}\label{prop:unik}
Let $\Ome$ be an unbounded domain in $\R^2$ such that the $r$ coordinate is
bounded on $\Ome$. Let $u$ and $v$ be two solutions of \eqref{eq:cmc1} on $\Ome$
which are continuous up to the boundary of $\Ome$ and such that $u=v$ along this
boundary. If the function $|v-u|$ is bounded on $\Ome$, then $u=v$.
\end{prop}

The proof is based on the same ideas as Theorem~2 in \cite{CoKu}
\begin{proof}
Let us define $\Ome_a=\{p\in\Ome\,|\ |p|<a\}$ and $C_a=\{p\in\Ome\,|\ |p|=a\}$.
let us define $w=v-u$ and $X=\frac{\cosh r\nabla v}{\sqrt{1+\cosh^2 r|\nabla v
|^2}}- \frac{\cosh r \nabla u}{\sqrt{1+\cosh^2 r|\nabla u |^2}}$. We denote by
$\vec\eta$ the outgoing normal to $\Ome_a$.

We then have
$$
\int_{C_a}w (\cosh r X)\cdot \vec\eta=\int_{\partial\Ome_a}w (\cosh r X)\cdot
\vec\eta=\int_{\Ome_a}\nabla
w\cdot(\cosh r X)
$$
By Lemma~1 in \cite{CoKu},
\begin{align*}
\nabla w\cdot (\cosh r X)&=(\cosh r \nabla v-\cosh r\nabla
u)\cdot\Big(\frac{\cosh r\nabla
v}{\sqrt{1+\cosh^2 r|\nabla v |^2}}- \frac{\cosh r \nabla u}{\sqrt{1+\cosh^2
r|\nabla u |^2}}\Big)\\
&\ge \left|\frac{\cosh r\nabla
v}{\sqrt{1+\cosh^2 r|\nabla v |^2}}- \frac{\cosh r \nabla u}{\sqrt{1+\cosh^2
r|\nabla u |^2}}\right|^2\\
&\ge |X|^2
\end{align*}

Let $a_0>0$ be such that $\Ome_{a_0}\neq\emptyset$ and denote
$\mu=\int_{\Ome_{a_0}}\nabla w\cdot(\cosh r X)$; $\mu>0$ if $u\neq v$. Since $w$
is bounded by a constant $M$ and $r$ is bounded on $\Ome$, $\cosh r$ is bounded
also by $M$. We then have:
$$
\mu+\int_{\Ome_a\setminus \Ome_{a_0}}|X|^2\le M^2\int_{C_a}|X|
$$
Let us define $I(a)=\int_{C_a}|X|$ and $l_a$ be the length of $C_a$; we
remark that $l_a\le 2\pi a$. We have
$$
I^2(a)=\Big(\int_{C_a}|X|\Big)^2\le l_a\int_{C_a}|X|^2
$$
Then 
\begin{equation}\label{eq:equadiff}
\mu+\int_{a_0}^a\frac{I^2(t)}{2\pi t}\dd t\le M^2 I(a)
\end{equation}

Let $\zeta$ be the function defined on $[a_0,a_0\exp(4\pi M^4/\mu))$ by the
following equation:
$$
-\frac{1}{\zeta}+\frac{2M^2}{\mu}=\frac{1}{2\pi M^2}\ln\frac{a}{a_0}
$$
$\zeta$ satisfies $\zeta(a_0)=\mu/(2M^2)$ and $M^2\zeta'=\zeta^2/(2\pi a)$.
Equation~\eqref{eq:equadiff} then implies that $I(a)\ge \zeta(a)$. But this is
impossible since $\zeta(a)$ converge to $+\infty$ when $a\rightarrow
a_0\exp(4\pi M^4/\mu)$. Then $u$ and $v$ are equal.
\end{proof}

\section{The monotonicity of the Alexandrov function}
\label{sec:monotonous}
This section is entirely devoted to the proof of Proposition~\ref{prop:decrease}. This
will finish the proof of our main result (Theorem~\ref{th:main}).


\subsection{The geometric configuration}

Let us consider a properly embedded annular top end $A$ with cmc $H_0$ which is
cylindrically bounded. Let $\Gamma=(\gamma_t)_{t\in\R}$ be a translation
foliation of $\H^2$ by geodesic lines. We assume that $\alpha_\Gamma$ is not
decreasing. 

By considering only $A\cap\{z\ge z_0\}$ for some $z_0>0$ large, we can assume
that $\alpha_\Gamma$ is increasing. Because of Proposition~\ref{prop:asymp}, we
can also assume that any horizontal section $A\cap\{z=z'\}$ is composed of one
curve with curvature strictly larger than $1$.

\begin{figure}[h]
\begin{center}
\resizebox{0.5\linewidth}{!}{\input{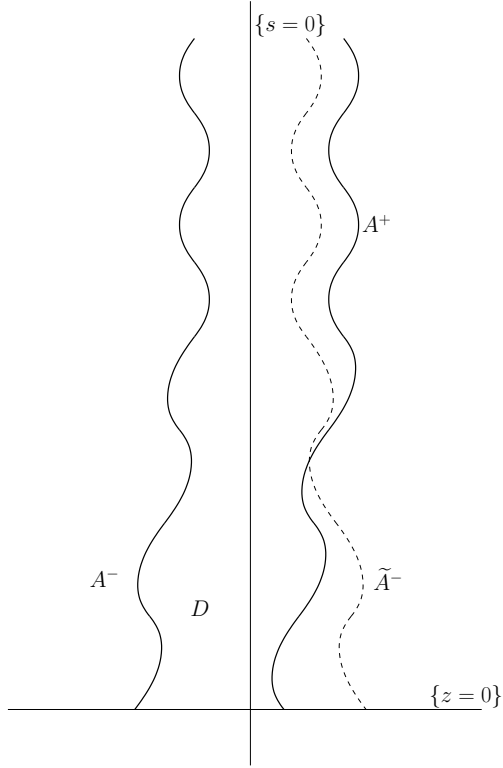}}
\caption{The annular end $A$}\label{fig:geoconf}
\end{center}
\end{figure}

Using the model introduced in Section~\ref{sec:horigraph} for $\H^2$, we can
assume that the foliation $\Gamma$ is given by $\gamma_t=\{s=t\}$. Moreover, by
changing the origin of the $s$ variable, we can assume that $\alpha_\Gamma(0)<0$
and $\lim_{+\infty}\alpha_\Gamma>0$. We also can assume that the intersection
$A\cap\{s=0\}$ is transverse. We define $A^+=A\cap\{s>0\}$ and
$A^-=A\cap\{s<0\}$ (see Figure~\ref{fig:geoconf}). 

The idea of the proof of Proposition~\ref{prop:decrease} is to obtain a control
of the flux of $A$ along $\partial A^+$ in the direction $\partial_s$. This idea
comes from the Positive Flux Lemma in \cite{KoKuMeSo}.

\subsection{A Dirichlet problem}\label{sec:dirichlet}

The annulus $A$ bounds a cylindricaly bounded domain $D$ in $\H^2\times\R_+^*$.
Let $\Ome$ denote the domain $\{(r,z)\in\R\times\R_+^*\,|\ (0,r,z)\in
D\}$. 

Let us denote by $\widetilde A^-$ the symmetric of $A^-$ by $\{s=0\}$. We then
define on $\Ome$ the function $f$ by $f(r,z)=\inf\{s\in\R_+\,|\ (s,r,z)\in
A^+\cup\widetilde A^-\}$. Since $A$ is cylindrically bounded, $f$ is uniformly
bounded on $\Ome$. Since the curvature of the curve $A\cap\{z=z'\}$ is larger
than $1$, $f$ extend continuously to $\partial \Ome$ and this boundary value is
$0$ along $\partial \Ome \cap \{z>0\}$.

Since $\alpha_\Gamma(0)<0$, the reflection procedure described in Section
\ref{sec:alexfunc} implies that $A^+\cap\{z=0\}$ is the horizontal Killing graph
of $f$ over $\partial\Ome\cap\{z=0\}$.

Besides, if $p\in\partial\Ome\cap\{z>0\}$ and the tangent space to $A$ at
$(0,p)$ is not normal to $\{s=0\}$, for $q\in\Ome$
$$
\frac{f(q)}{d(q,\partial\Ome)} \text{ is bounded in a neighborhood of }p
$$
where $d(q,\partial\Ome)$ is the Euclidean distance from $q$ to $\partial\Ome$.

The aim of this section is to solve a Dirichlet problem for \eqref{eq:cmc1}: we
prove the following result.
\begin{lem}\label{lem:dirichlet}
There exists a unique nonnegative solution $u$ on $\Ome$ of
Equation~\eqref{eq:cmc1} which is continuous up to the boundary with boundary
value $f$ and such that $u\le f$ on $\Ome$.
\end{lem}

\begin{proof}
The uniqueness comes from Proposition~\ref{prop:unik} since the $r$ coordinate is bounded in $\Ome$ ($D$ is cylindrically bounded).

For the existence part, we use the Perron method to solve the Dirichlet
problem. Let us recall the framework of the Perron method.

Let $v$ be a continuous function on $\overline\Ome$, $v$ is called a subsolution
for our problem if $v\le f$ and if, for any compact subdomain $U\subset \Ome$ and
any solution $h$ of \eqref{eq:cmc1} with $v\le h$ on the boundary of $U$, we have
$v\le h$ on $U$. If $S$ denote the set of all subsolutions, we define our
solution by the following formula:
$$
u(q)=\sup_{v\in S} v(q)
$$
We notice that $S$ is non empty since the function $0$ is a subsolution; thus
$u\ge 0$.
Moreover if $v$ and $w$ are subsolutions, the continuous function $\max(v,w)$ is
also a subsolution. It is also clear that $u\le f$ but it is not clear that $u$
is a solution to our problem.

Since the $r$ coordinate is bounded in $\Ome$, Proposition~\ref{prop:exist} implies
there is a $R>0$ such that, for any disk $\Delta\subset \Ome$ of radius less
than $R$, the Dirichlet problem can be solved on $\Delta$ for
Equation~\eqref{eq:cmc1}. Thus for any such disk $\Delta$ and subsolution $v$,
we can define the continuous function $M_\Delta(v)$ on $\Ome$ as $M_\Delta(v)=v$
outside $\Delta$ and $M_\Delta(v)$ is equal to the solution of \eqref{eq:cmc1}
in $\Delta$ with $v$ as boundary value. Since $v$ is
a subsolution, $v\le M_\Delta(v)$. The graph of $v$ is below $A^+$ and
$\widetilde A^-$ so, by the maximum principle for cmc $H_0$ surfaces, the graph
of $M_\Delta(v)$ is below $A^+$ and $\widetilde A^-$. This implies that
$M_\Delta(v)$ is a subsolution. 

\begin{claim}\label{cl:perron1}
The function $u$ is a solution of \eqref{eq:cmc1} in $\Ome$.
\end{claim}

\begin{proof}[Proof of Claim \ref{cl:perron1}]
Let us consider $p\in\Ome$ and $\Delta$ a disk in $\Ome$ centered at $p$ with radius
less than $R$. Let $(v_n)$ be a sequence of subsolutions such that
$v_n(p)\rightarrow u(p)$. By considering $\max(0,v_n)$ we can assume $v_n\ge 0$.
We have $M_\Delta(v_n)$ is also a sequence of subsolutions with
$M_\Delta(v_n)(p)\rightarrow u(p)$. On $\Delta$, $M_\Delta(v_n)$ is a bounded
sequence of solutions of \eqref{eq:cmc1}. So by considering a subsequence if
necessary, we can assume that it converges to a solution $\bar v$ on $\Delta$
with $u\ge \bar v$ and $\bar v(p)=u(p)$. Let us prove that $u=\bar v$ on
$\Delta$, so $u$ will be a solution of \eqref{eq:cmc1}.

If it is not the case there is a point $q\in \Delta$ where $u(q) >\bar v(q)$. So
there is a subsolution $w$ such that $w(q)>\bar v(q)$. So let us consider the
sequence of subsolutions $M_\Delta(\max(w,v_n))$. We have
$M_\Delta(\max(w,v_n))\ge M_\Delta(v_n)$ and $M_\Delta(\max(w,v_n))\ge w$.
Moreover on $\Delta$, it is a sequence of solutions of \eqref{eq:cmc1}; so
considering a subsequence, it converges to a solution $\bar w$ of
\eqref{eq:cmc1} with $\bar w\ge \bar v$ and $\bar w(q)\ge w(q)>\bar v(q)$. But
since $\bar w(p)=\bar v(p)$, the maximum principle gives $\bar w=\bar v$ on
$\Delta$ which contradicts $\bar w(q)\ge w(q)>\bar v(q)$. The claim is proved.
\end{proof}

Since $0\le u\le f$, the function $u$ is continuous up to the part of
$\partial\Ome$ not in $\{z=0\}$ and $u=0=f$ there. For
$\partial\Ome\cap\{z=0\}$, we need to construct some barriers for the problem.

\begin{claim}\label{cl:barrier}
Let $r_0$, $M>0$ and $\eps>0$ be real numbers. There exist a neighborhood $V$ of
$(r_0,0)$ in $\R\times\R_+$ as small as we want and a solution $h$ of
\eqref{eq:min} in $V$ which is continuous up to $\partial V$ such that $0\le
h\le M$ on $V$, $h=0$ on $\partial V\cap (\R\times\R_+^*)$ and
$h(r_0,0)=M-\eps$.
\end{claim}

\begin{proof}[Proof of Claim \ref{cl:barrier}]
Let us consider $r_0\in\R$, $M>0$, $\eps>0$ as in the claim and $R>0$ such that
the Dirichlet problem for the minimal surface equation \eqref{eq:min} can be
solved on the disk $\Delta$ centered at $(r_0,0)$ and radius $R$
(Proposition~\ref{prop:exist}). On
$\partial \Delta$, let $\phi_n$ be a continuous function such that 
\begin{itemize}
\item $0\le \phi_n\le 2M$.
\item $0=\phi_n$ on $\partial \Delta\cap\{z>1/n\}$ and $\phi_n\le M$ on $\{z\ge
0\}$.
\item $\phi_n(r,z)=2M-\phi_n(r,-z)$.
\end{itemize}
Moreover, we assume $\phi_n\ge\phi_{n+1}$ on $\partial \Delta\cap\{z\ge0\}$. Let
$h_n$ be the solution of \eqref{eq:min} on $\Delta$ such that
$h_n=\phi_n$ on $\partial \Delta$. By uniqueness of the solution and the maximum
principle, we have
\begin{itemize}
\item $0\le h_n\le 2M$.
\item $h_n(r,z)=2M- h_n(r,-z)$ so $h_n(r,0)=M$.
\item $h_n$ is decreasing and $h_n\le M$ on $\Delta\cap\{z\ge 0\}$
\end{itemize}

Since the sequence is bounded, it converges to a solution $\bar h$ of
\eqref{eq:min} on $\Delta$. Because of the monotonicity, $\bar h$ is continuous up
to the boundary except at the points $(r_0+R,0)$ and $(r_0-R,0)$. $\bar h$ is
equal to $0$ on $\partial \Delta\cap\{z>0\}$ and to $M$ on $\Delta\cap\{z=0\}$. So by
continuity there is an $\eta>0$ such that $\bar h(r_0,\eta)=M-\eps$. So if we
consider the restriction of $\bar h$ to $\Delta\cap\{z\ge \eta\}$ we have constructed
a neighborhood $V$ of $(r_0,\eta)$ in $\{z\ge \eta\}$ and a solution $\bar h$ of
\eqref{eq:min} on $V$ which is continuous up to the boundary such
that $0\le \bar h\le M$ on $V$, $\bar h=0$ on $\partial V\cap \{z>\eta\}$ and
$\bar h(r_0,\eta)=M-\eps$. We notice that by choosing $R$ small, we can assume
$V$ as small as we want.
\end{proof}

With these barriers we can finish the proof of Lemma~\ref{lem:dirichlet}.

\begin{claim}\label{cl:continuity}
The function $u$ is continuous up to the boundary of $\Ome$ and takes the value
$f$ on it.
\end{claim}

\begin{proof}[Proof of Claim \ref{cl:continuity}]
The problem is only on $\partial\Ome\cap\{z=0\}$ minus its end points; so take a
point $p\in \partial\Ome\cap\{z=0\}$.
Let us consider $\eps>0$ and $I$ a segment in $\{z=0\}$ containing $p$ such that
$f\ge f(p)-\eps$ on $I$. Now from our construction of barriers, we know that
there exist a neighborhood $V$ of $p$ in $\Ome$ such that $\overline V\cap
\{z=0\}\subset I$ and a solution $h$ of \eqref{eq:min} on $V$
continuous up to the boundary and such that $h=0$ on $\partial V\cap \Ome$,
$h\le f(p)-\eps$ on $\partial V\cap \{z=0\}$ and $h(p)=f(p)-2\eps$. Let us
extend the definition of $h$ by $0$ to the whole $\Ome$. By the maximum
principle, $h$ is a subsolution for our problem, so $u\ge h$. This implies that
$\liminf_p u\ge f(p)-2\eps$. Since $u\le f$ on $\Ome$ we have $\limsup_p u\le
f(p)$. Then $u$ is continuous at $p$ and takes the value $f(p)$.
\end{proof}
\end{proof}


\subsection{The asymptotic behaviour of $u$}

We know from Proposition~\ref{prop:asymp} that the annular end $A$ is asymptotic
for large $z$ to a Delaunay surface. In this subsection, we will see that this
asymptotic behavior passes to the function $u$.

Let $(z_n)$ be a sequence such that $z_n\nearrow +\infty$ and $t_{-z_n}(A)$
converges to a Delaunay surface $\boD_\tau$, let us denote by $G$ the cylindrically
bounded domain whose boundary is $\boD_\tau$. We notice that by our normalization
of $A$, the axis of $\boD_\tau$ is $\{r=0,s=\lim_{+\infty}\alpha_\Gamma\}$.

Let us also denote by $t_{-z_n}$ the translation by $-z_n$ in the $(r,z)$ plane.
Because of the asymptotic behavior of $A$, the sequence of domains
$t_{-z_n}(\Ome)$ converges to the domain $\Ome_0$ defined by
$\Ome_0=\{(r,z)\in\R^2\,|\ (0,r,z)\in G\}$ (the convergence is smooth on any
compact). Let us define $\boD_\tau^-=\boD_\tau\cap \{s<0\}$ and
$\widetilde{\boD}_\tau^-$ the symmetric of $\boD_\tau^-$ by $\{s=0\}$. $\widetilde
\boD_\tau^-$ is a horizontal Killing graph of a function $f_0$ on $\Ome_0$.
Actually, $f_0=\lim f\circ t_{z_n}$.

Let us consider the sequence $u\circ t_{z_n}\le f\circ t_{z_n}$ on
$t_{-z_n}(\Ome)$, it is a uniformly bounded sequence. So if we consider a
subsequence, we can assume that $u\circ t_{z_n}$ converges to a solution $u_0$ of
\eqref{eq:cmc1} on $\Ome_0$. Moreover we have $0\le u_0\le f_0$. This implies
that $u_0$ is continuous up to $\partial\Ome_0$ and takes value $0$ there. We
then have $u_0$ and $f_0$ two solutions of \eqref{eq:cmc1} on $\Ome_0$ with the
same vanishing boundary value; so, by Proposition~\ref{prop:unik}, $u_0=f_0$.

The uniqueness of the possible limit implies that the whole sequence $u\circ
t_{z_n}$ converges to $f_0$.


\subsection{Computation of fluxes}

The idea of this section is to compute the flux of $A$ along the boundary of
$A^+$ in the direction of $\partial_s$ and find a contradiction which will prove
Proposition~\ref{prop:decrease}.

Let $(z_n)$ be an increasing sequence such that
$z_n\rightarrow +\infty$ and $\Ome_0$ be the associated limit domain. This
domain is either a strip if $\tau=2H_0-\sqrt{4H_0^2-1}$ or a periodic domain
composed of successive "bubbles" if $0<\tau<2H_0-\sqrt{4H_0^2-1}$. By adding a
constant to $(z_n)$, we assume that $\{z=0\}$ is a line of symmetry of
$\Ome_0$. If $u_0=\lim u\circ t_{z_n}$, we get that $u_0(r,z)=u_0(r,-z)$ and
$\boD_\tau$ is symmetric with respect to $\{z=0\}$.

The boundary of $A_n^+=A^+\cap \{0<z<z_n\}$ is composed of four smooth arcs:
$\gamma^1= A^+\cap\{z=0\}$, $\gamma_n^2= A^+\cap \{z=z_n\}$ and
$\gamma_n^3=A\cap\{s=0\}\cap\{0\le z\le z_n\}$ ($\gamma_n^3$ is actually
composed of two arcs). The flux of $A$ along $\partial A_n^+$ in the direction
of $\partial_s$ is equal to $0$ since $\partial A_n^+$ is homologically trivial.
The idea is to use the graph of $u$ as a barrier for the computation of
$F_{\partial A_n^+}(\partial_s)$ to prove that it can not vanish for large $n$.

Let us denote by $\Ome_{z_n}$ the subdomain $\Ome\cap\{0<z<z_n\}$. In order to
compute the flux, we need a surface $Q$ bounded by $\partial
A_n^+$: we define $Q$ as the union of $G\cap\{s\ge 0,z=0\}$, $G\cap\{s\ge
0,z=z_n\}$ and $\{0\}\times\Ome_{z_n}$. The term $(\partial_s,\vec n_Q)$ is zero
along the first two parts and is equal to $-\cosh r$ along the third part so
the flux of $A$ along $\partial A_n^+$ is equal to 
\begin{equation}\label{eq:flux1}
0=F_{\partial A_n^+}(\partial_s)=\int_{\partial
A_n^+}(\vec\nu,\partial_s)+2H_0\int_{\Ome_{z_n}}\cosh r \dd r \dd z
\end{equation}

On an other hand, we have from \eqref{eq:cmc1}
\begin{equation}\label{eq:flux2}
0=\int_{\partial\Ome_{z_n}}(\frac{\cosh^2 r \nabla u}{\sqrt{1+\cosh^2 r |\nabla
u|^2}},\vec\eta) \dd s+\int_{\Ome_{z_n}}2H_0\cosh r\dd r\dd z
\end{equation}
with $\vec\eta$ the outgoing unit normal to $\Ome_{z_n}$.
We notice that, even if we do not know that $u$ is smooth up to the boundary of
$\Ome$, the first integral is well defined since the vector field $\frac{\cosh^2
r \nabla u}{\sqrt{1+\cosh^2 r |\nabla u|^2}}$ is bounded and
Equation~\eqref{eq:cmc1} is satisfied.

Thus Equations~\eqref{eq:flux1} and \eqref{eq:flux2} give
\begin{equation}\label{eq:flux3}
0=\int_{\partial\Ome_{z_n}}(\frac{\cosh^2 r \nabla u}{\sqrt{1+\cosh^2 r |\nabla
u|^2}},\vec\eta) \dd s -\int_{\partial A_n^+}(\vec\nu,\partial_s)
\end{equation}

In order to compare the terms in the above equality, we need to study the
regularity of $u$ at the boundary of $\Ome$.

Let us define $c^1=\partial\Ome_{z_n}\cap \{z=0\}$, $c_n^2=\partial
\Ome_{z_n}\cap \{z=z_n\}$ and $c_n^3=\partial \Ome_{z_n}\setminus(c^1\cup
c_n^2)$. We notice that $\gamma_n^3=\{0\}\times c_n^3$, so the same notation can
be used to denote the two curves.

\begin{claim}
The function $u$ is $C^2$ up to the boundary at each point of $c^1$ (except its
end-points) and each point in $c_n^3$ where $A$ is not normal to $\{s=0\}$.
\end{claim}

\begin{proof}
Let $p$ be a point $c_n^3$ where $A$ is not normal to $\{s=0\}$. As written at
the beginning of subsection~\ref{sec:dirichlet}, there is a neighborhood of $p$
such that $f/d(\cdot,\partial\Ome)$ is
bounded. Since $0\le u \le f$, Proposition~\ref{prop:estim} gives a uniform
upper-bound for $|\nabla u|$ in a neighborhood of $p$. Since $u$ satisfies
Equation~\eqref{eq:cmc1} and the boundary data are smooth, elliptic regularity
theory implies that $u$ is $C^2$ up to the boundary near $p$ (see for example
\cite{LaUr}, Theorem~4.6.1 gives $C^{1,\alpha}$ regularity and Theorem~4.6.3
gives $C^2$ regularity).

For the regularity at a point in $c^1$, the argument is the same but we have to
prove the uniform upper bound for the gradient. We notice that, near $\gamma^1$,
$A^+$ is the graph of $f$. So, near $c^1$, $f$ is a smooth function. Since $u\le
f$ this give a barrier from above with bounded gradient for $u$. The problem is
to obtain a barrier from below.

The curves $\gamma^1$ and $\{0\}\times c^1$ bounded a convex domain in
$\H^2\times \{0\}$ whose boundary is thus composed of an arc of curvature larger
than $1$ and a geodesic arc.
Let us denote by $U$ this domain viewed in $\H^2\times \{0\}$. Let $q$ be the
middle of $\{0\}\times c^1$. Let $\eps$ be positive and $\Gamma_\eps$ be the
Jordan arc in $\H^2\times \R$ composed by $\gamma^1$ and two geodesic arc
joining the end points of $\gamma^1$ to $t_\eps(q)$. When $\eps$ is small, the
two geodesic arcs are included in $\{0\}\times \Ome$. Let $\Sigma_\eps$ be the
solution of the Plateau problem in $\H^2\times\R$ with $\Gamma_\eps$ as
boundary. Since $\Gamma_\eps $ is a vertical graph above the boundary of the
convex domain $U$, $\Sigma_\eps$ is unique and is a vertical graph above $U$.
Moreover, this graph is smooth up to the boundary (barriers from above and below
can be easily found). $\Sigma_\eps$ is included in $\{z\ge 0\}$ and can not be
tangent to $\{z=0\}$ by the maximum principle. The translate
$t_{-a}(\Gamma_\eps)$ for $a>0$ never meets the graph of $u$. So by the maximum
principle, $t_{-a}(\Sigma_\eps)$ never meets the graph of $u$ and then
$\Sigma_\eps$ is inside $\{(s,r,z)\in\R\times \Ome\, |\ s\le u(r,z)\}$. Since
$\Sigma_\eps$ is not tangent to $\{z=0\}$ along $\gamma^1$, $\Sigma_\eps$ is a
good barrier from below for $u$.  

We then get a uniform bound for the gradient of $u$ near any points of $c^1$
(Proposition~\ref{prop:estim}); this gives us the $C^2$ regularity up to the
boundary.
\end{proof}

Using the regularity of the function $u$, we prove the following statement.

\begin{claim}\label{cl:flux1}
We have:
$$
\int_{c^1}(\frac{\cosh^2 r \nabla u}{\sqrt{1+\cosh^2 r |\nabla
u|^2}},\vec\eta)  -\int_{\gamma^1}(\vec\nu,\partial_s)>0
$$
and
$$
\int_{c_n^3}(\frac{\cosh^2 r \nabla u}{\sqrt{1+\cosh^2 r |\nabla
u|^2}},\vec\eta) -\int_{\gamma_n^3}(\vec\nu,\partial_s)>0
$$
\end{claim}

\begin{proof}
The curve $\gamma^1$ is the graph of the function $f$ over $c^1$ and, in fact,
close to $\gamma^1$, $A^+$ is the graph of $f$ so we have 
$$
\int_{\gamma^1}(\vec\nu,\partial_s)=\int_{c^1}(\frac{\cosh^2 r \nabla
f}{\sqrt{1+\cosh^2 r |\nabla f|^2}},\vec\eta) 
$$

Now since $u\le f$ and by the maximum principle $u$ and $f$ can not have the
same gradient on $c^1$ this implies that along $c^1$ we have
$$
(\frac{\cosh^2 r \nabla u}{\sqrt{1+\cosh^2 r |\nabla
u|^2}},\vec\eta)>(\frac{\cosh^2 r \nabla f}{\sqrt{1+\cosh^2 r |\nabla
f|^2}},\vec\eta) 
$$
This give the first inequality.

For a point $p$ in $\gamma_n^3$,
if $A$ is normal to $\{s=0\}$ $(\vec\nu,\partial_s)=-\cosh r$. Since
$\frac{\cosh^2 r \nabla u}{\sqrt{1+\cosh^2 r |\nabla u|^2}}$ has a norm less
than $\cosh r$ everywhere, we get a large inequality at $p$ between the to
integrand.

If $A$ is not normal to $\{s=0\}$ at $p$, we remark that the term
$(\vec\nu,\partial_s)$ give the same result if $\vec \nu$ is the conormal to
$A^+$ or $\widetilde A^-$. Besides the graph of $u$ is regular up to the
boundary and is below $A^+$ and $\widetilde A^-$. So, with $\vec\nu_u$ the
conormal to the graph of $u$ and $\vec\nu_A$ the one for $A^+$, the maximum
principle implies that $(\vec\nu_u,\partial_s)>(\vec \nu_A,\partial_s)$ at $p$.
After integration we get the second inequality since there is always point where
$A$ is not normal to $\{s=0\}$ on $\gamma_n^3$.
\end{proof}

We have the following limits for the integrals along $\gamma_n^2$ and $c_n^2$.

\begin{claim}\label{cl:flux2}
We have the following limits:
$$
\lim_{n\rightarrow \infty}\int_{\gamma_n^2}(\vec\nu,\partial_s)=0 
$$
and
$$
\lim_{n\rightarrow \infty}\int_{c_n^2}(\frac{\cosh^2 r \nabla u}{\sqrt{1+\cosh^2
r |\nabla u|^2}},\vec\eta)=0
$$
\end{claim}

\begin{proof}
Since $t_{-z_n}(A)$ converges to $\boD_\tau$, the first limit is equal to the
integral of $(\vec\nu,\partial_s)$ along the curve $\boD_\tau\cap\{s\ge 0,z=0\}$.
By our choice of $(z_n)$, the surface $\boD_\tau$ is symmetric with respect to
$\{z=0\}$. The conormal is then equal to $\partial_z$ and the scalar product
vanishes. The limit is then $0$.

For the second limit, we know that $u\circ t_{z_n}$ converges to $u_0$ this
implies that the limit of the integral is equal to
$$
\int_{\Ome_0\cap\{z=0\}}(\frac{\cosh^2 r \nabla u_0}{\sqrt{1+\cosh^2 r |\nabla
u_0|^2}},\partial_z)
$$
We notice that \textit{a priori} the convergence of $u\circ t_{z_n}$ is smooth
only on compact subdomains of $\Ome_0$ but since the integrand is uniformly
bounded it is sufficient to take the limit of the integral. Now, we have
$u_0(r,z)=u_0(r_0,-z)$ so $\nabla u_0$ is normal to $\partial_z$ on $\Ome_0\cap
\{z=0\}$ and the limit integral vanishes.
\end{proof}

Now using Equation~\eqref{eq:flux3} and Claims~\ref{cl:flux1} and
\ref{cl:flux2}, we get our contradiction which finishes the proof of
Proposition~\ref{prop:decrease}
\begin{align*}
0&=\lim \Big(\int_{\partial\Ome_{z_n}}(\frac{\cosh^2 r \nabla u}{\sqrt{1+
\cosh^2 r |\nabla u|^2}},\vec\eta) \dd s -\int_{\partial A_n^+}(\vec\nu,
\partial_s)\Big)\\
&\ge \int_{c^1}(\frac{\cosh^2 r \nabla u}{\sqrt{1+\cosh^2 r |\nabla
u|^2}},\vec\eta) \dd s -\int_{\gamma^1}(\vec\nu,\partial_s)\\
&>0
\end{align*}

\bibliographystyle{plain}
\bibliography{../reference}

\end{document}